\documentclass{IEEEtran}
\usepackage{textcomp}
\usepackage{epsfig}
\usepackage{color}
\usepackage{amssymb,latexsym,amsfonts,amsmath,mathrsfs}
\usepackage{graphicx}
\usepackage{cite}
\usepackage{algorithmic} 
\usepackage{bbm}
\newtheorem{thm}{Theorem}
\newtheorem{lem}{Lemma}

\newtheorem{corollary}{Corollary}

\newtheorem{assum}{Assumption}
\newtheorem{rem}{Remark}

\newtheorem{exmp}{Example}
\def\BibTeX{{\rm B\kern-.05em{\sc i\kern-.025em b}\kern-.08em
		T\kern-.1667em\lower.7ex\hbox{E}\kern-.125emX}}
\IEEEoverridecommandlockouts
\usepackage{epsfig}
\usepackage{color}
\usepackage{amssymb,latexsym,amsfonts,amsmath,mathrsfs}
\usepackage{graphicx}
\usepackage{cite}
\usepackage{algorithmic} 
\usepackage{bbm}

\usepackage{comment}

\usepackage[linesnumbered,ruled]{algorithm2e}
\allowdisplaybreaks[4]

\begin{document}
\allowdisplaybreaks[4]
	
\title{Distributed Stochastic Zeroth-Order Optimization with Compressed Communication}
\author{Youqing Hua,  Shuai Liu, \IEEEmembership{Member, IEEE}, Yiguang Hong, \IEEEmembership{Fellow, IEEE}, and Wei Ren, \IEEEmembership{Fellow, IEEE}}

\maketitle

\begin{abstract}
The dual challenges of high communication costs and gradient inaccessibility--common in privacy-sensitive systems or black-box environments--drive our work on communication-constrained, gradient-free distributed optimization. We propose a compressed distributed stochastic zeroth-order algorithm (Com-DSZO), which requires only two function evaluations per iteration and incorporates general compression operators. Rigorous analysis establishes a sublinear convergence rate for both smooth and nonsmooth objectives, explicitly characterizing the trade-off between compression and convergence. Furthermore, we develop a variance-reduced variant (VR-Com-DSZO) under stochastic mini-batch feedback. The effectiveness of the proposed algorithms is demonstrated through numerical experiments.
\end{abstract}


\section{Introduction}\label{sec:introduction}
Over the past decade, the extensive application of Distributed Convex Optimization (DCO) in Multi-Agent Systems (MASs)—particularly in smart grids \cite{smartgrid}, machine learning \cite{machinelearning}, and sensor networks \cite{gt}—has significantly advanced both practical implementations and theoretical research. In these systems, each distributed agent is endowed with a private loss function. The primary objective of the MAS is to minimize the average of these local loss functions through local interactions.\par
To date, several distributed algorithms have been developed \cite{sundharram_distributed_2010,DSGT,DPGD,b4,ADMM,b38,gt,b37}, including distributed projected sub-gradient methods \cite{sundharram_distributed_2010,DPGD,b4}, distributed ADMM \cite{ADMM}, and distributed gradient-tracking algorithms \cite{gt}. Unlike the deterministic methods proposed in \cite{DPGD,b4,ADMM,b38,b37}, the studies in \cite{sundharram_distributed_2010,gt,DSGT} focus on distributed stochastic optimization, where each local loss function is defined as the expectation of a random function. This class of problems has garnered considerable interest due to its wide applicability. A common paradigm for solving such problems is to combine gradient-based algorithms with consensus protocols, such as Distributed Stochastic Gradient Descent (DSGD) \cite{sundharram_distributed_2010}. However, in fields such as biochemistry, geoscience, and aircraft circling, accessing gradients may pose significant challenges due to data privacy concerns, black-box constraints, or the high computational costs involved.\par 
The Zeroth-Order (ZO) method, which estimates gradients via function evaluations, offers a practical alternative to first-order algorithms. A seminal work by Flaxman et al. introduced a one-point gradient estimator \cite{flexman-onepoint}, followed by the development of two-point sampling ZO algorithms \cite{nesterov2017, b14} and extensions to stochastic constrained optimization \cite{duchi}. More recently, ZO algorithms have been extended to distributed settings \cite{DmatchC, nonsmooth, noncvx-dso-zo, online1, online2, online3, b22, b23, zero-ordernonconvex}. In deterministic DCO, researchers have analyzed the convergence of distributed ZO gradient descent algorithms \cite{DmatchC} and addressed non-smooth problems \cite{nonsmooth}. In stochastic distributed optimization, primal-dual algorithms for nonconvex objectives have been proposed \cite{noncvx-dso-zo}. Distributed ZO optimization has also been investigated in online settings with time-varying loss functions \cite{online1, online2, online3}. Despite these advances, the need for high-precision iterative communication between agents remains a significant bottleneck in distributed optimization.\par 
To address this issue, considerable attention has been devoted to developing effective communication compression schemes for distributed algorithms. A compression operator \( \mathcal{C}(\cdot)\) is typically employed to produce a bit-coded output that substantially reduces communication overhead relative to the original input.
Significant progress has been made in distributed first-order optimization \cite{b25,b26,unbiased1,Choco-gossip,onoine-choco,b29,b30,b33}. Early work introduced DSGD algorithms with sparse gradients \cite{b25} and randomized quantized gradients \cite{b26}, although these approaches were initially limited to specific classes of compressors. Subsequent research investigated contractive stochastic compressors under the assumption of unbiasedness \cite{unbiased1}, which implicitly excludes many widely used deterministic compressors. In contrast, the Choco-Gossip algorithm \cite{Choco-gossip} was developed to accommodate biased yet contractive stochastic compressors and was later extended to online DCO scenarios \cite{onoine-choco}. Further advancements expanded the scope to a broader class of stochastic compressors \cite{b29,b30,b33} under static parameters for unconstrained optimization. Specifically, \cite{b29} and \cite{b30} introduced gradient-tracking-based methods that feature explicit trade-offs between compression and convergence for $\delta$-contracting compressors, which necessitate $L$-smoothness and $\mu$-strong convexity conditions. In contrast, \cite{b33} provided convergence guarantees but lacks a detailed rate analysis. Notably, current literature exhibits a critical gap in designing distributed constrained optimization algorithms with adaptive step-size policies for general compressors.\par
Although distributed first-order optimization with compression has been extensively studied, research on distributed ZO algorithms with compression remains limited and underdeveloped. The most relevant work \cite{unbiased1} proposed a distributed online ZO algorithm based on one-point gradient estimation, but it is restricted to unbiased compressors and smooth functions. The dual challenges of bias and variance in stochastic gradient estimation underscore the need for further advances in communication-efficient distributed ZO optimization.\par 
This paper addresses these challenges by developing a communication-efficient algorithm for stochastic distributed ZO optimization. We consider a broad class of compressors and analyze the convergence rates of the proposed algorithm for both smooth and nonsmooth objectives. The main contributions are summarized as follows:
\begin{enumerate}
		\item We propose Com-DSZO, a communication-efficient distributed zeroth-order optimization framework for stochastic constrained optimization in smooth/non-smooth settings. The method requires merely two function evaluations per iteration while eliminating gradient dependence, and demonstrates universal compatibility with compression operators, $\!$including the biased/unbiased compression schemes investigated in \cite{unbiased1,Choco-gossip,onoine-choco}.\par 
		\item We provide a detailed convergence analysis that explicitly quantifies the impact of compression resolution and problem dimensionality. In particular, for compressors with bounded relative error, Com-DSZO achieves a convergence rate of $\mathcal{O}\left( \frac{d^2\sqrt{T+c}}{( \psi r \omega)^2} \right)$, where $d$ denotes the problem dimension, $\omega \in (0,1]$ is the compression factor, and $r$ and $\psi$ satisfy $r \psi \le 1$. When $\omega \to 1$, corresponding to the uncompressed case, Com-DSZO recovers the convergence rate of DSZO \cite{DmatchC}.
		\item A variance-reduced variant, VR-Com-DSZO, is developed under stochastic mini-batch feedback to reduce variance in stochastic gradient estimation. Under the condition $b_1 b_2 = d \sqrt{d}$, VR-Com-DSZO attains the state-of-the-art convergence rate of DSGD methods \cite{sundharram_distributed_2010}, even in decentralized and communication-constrained environments.
\end{enumerate}\par 
\textit{Notation:} Let $\mathbb{R}^d$ denote the $d$-dimensional Euclidean space, and $\mathbb{S}^d$ the unit sphere, centered at the origin. The vectors \(\mathbf{1}_d\) and \(\mathbf{0}_d\), and the matrix \(\mathbf{I}_d\) denote the \(d\)-dimensional all-ones vector, the \(d\)-dimensional all-zeros vector, and the \(d \times d\) identity matrix, respectively. Subscripts are omitted when the context unambiguously identifies the dimension. $\|\!\cdot\!\|_2$ and $\|\!\cdot\!\|_F$ denote the Euclidean norm and the Frobenius norm of a matrix, respectively. The expectation with respect to the random vector $\xi$ is written as $\mathbb{E}_\xi[\cdot]$. The projection onto a convex set $\mathcal{C}$ is denoted by $\Pi_{\mathcal{C}}[\cdot]$ and defined as $\Pi_{\mathcal{C}}[z] := \arg\min_{x \in \mathcal{C}} \|x - z\|_2^2$. For symmetric matrices $\mathbf{M}$ and $\mathbf{N}$, $\mathbf{M} \succ \mathbf{N}$ indicates that $\mathbf{M} - \mathbf{N}$ is positive definite. The entry in the $i$-th row and $j$-th column of matrix $\mathbf{W}$ is denoted by $\mathbf{W}_{i,j}$, and $\text{null}(\mathbf{W})$ denotes its null space. For a differentiable function $f$, $\nabla f$ denotes its gradient. A function $f: \mathbb{R}^d \to \mathbb{R}$ is $L$-Lipschitz continuous if and only if $\left| f(x) - f(y) \right| \leq L \|x - y\|$ for all $x, y \in \mathbb{R}^d$. The notation $\beta(k) = \mathcal{O}(\alpha(k))$ means $\limsup_{k \to \infty} \frac{\beta(k)}{\alpha(k)} < +\infty$.\par
\section{Problem Formulation and Preliminaries}
\subsection{Stochastic Distributed Constrained Optimization}
This paper explores the stochastic constrained DCO problem of a MAS with $n$ agents. Each agent $i$ is associated with a decision vector $x_i\in \Omega \subseteq \mathbb{R}^d$ and an expectation-valued loss function $f_i\left( x \right)\triangleq\mathbb{E}_{\xi_i}\left[ {{F}_{i}}(x,\xi_i) \right]$, where $\xi_i$ is a random sample drawn from an unknown data distribution $\mathcal{P}_i$. The goal of the system is to cooperatively solve the following problem:
\begin{align}\label{1}
	\underset{x\in \Omega }{\mathop{\min }}\,f(x):=\frac{1}{n}\sum\nolimits_{i=1}^{n}{{{f}_{i}}\left( x \right)}\text{,}
\end{align}
where $x$ is the global decision vector and $f\left(x\right)$ is the global loss function. Moreover, for each agent $i$, there exists a random ZO oracle. Specifically, given a local decision vector $x_i$ and a random sample $\xi_i$, a private stochastic ZO oracle $F_i(x_i,\xi_i)$ can be observed only by agent $i$. No additional information, such as the gradient of \( F_i(x_i, \xi_i) \), is available.\par 
The information interaction among agents is modeled as an undirected graph $\mathcal{G}=(\mathcal{V},\mathcal{E})$, where $\mathcal{V}=\{1,2,\cdots,n\}$ stands for the nodes set, $\mathcal{E}\in \mathcal{V}\times \mathcal{V}$ is the edges set. The set of neighbors for agent $i$ is denoted by ${\mathcal{N}_i}$. The corresponding weighting matrix is $\mathbf{W}=[w_{i,j}]_{n\times n}$.
The standard assumptions on graph $\mathcal{G}$ are given below.\par
\begin{assum}\label{assum G}
	The undirected graph $\mathcal{G}$ is connected and balanced, and consequently the associated weight matrix $\mathbf{W}$ satisfies: \textit{a) symmetry}: $\mathbf{W}^\mathsf{T}=\mathbf{W}$; \textit{b) consensus property}: $\text{null}\left(\mathbf{I}-\mathbf{W}\right)=0$; \textit{c) spectral property:} $\mathbf{I}+\mathbf{W}\succ 0$.
\end{assum}\par 
\begin{rem}
	$\mathcal{G}$ is a balanced graph, implying that $\mathbf{W}$ is doubly stochastic; i.e., $\mathbf{W}\mathbf{1}=\mathbf{1}^\mathsf{T} \mathbf{W} =1$. Furthermore, based on Assumption \ref{assum G} and the Perron-Frobenius theorem, we can deduce that the largest singular value of $\mathbf{W}$ is equal to 1, while the remaining singular values are strictly less than 1. For analytical convenience, we denote the spectral gap of $\mathbf{W}$ by $\delta \in (0, 1]$ and let $\beta = \left\| \mathbf{I}-\mathbf{W} \right\|_2\in(0,2)$. Also, we have $\left\| \mathbf{W} - \mathbf{H} \right\|_2 = 1 - \delta$ and $\left\| \mathbf{H}\right\|_2 = 1$, where $\mathbf{H}=\frac{1}{n}\mathbf{1}\mathbf{1}^{\mathsf{T}}$.
\end{rem}\par 
For the optimization problem \eqref{1}, we assume the following.\par 
\begin{assum}\label{assum1-set}
	The set $\Omega$ is closed and convex, with bounds $\bar{R}>0$ and $\underline{R}>0$, ensuring $\underline{R}\le {{\left\| x \right\|}_{2}}\le \bar{R}$ for any $x\in \Omega $.\par
\end{assum}
\begin{assum}\label{assum2-f}
	For any \( i \in \mathcal{V} \), \( F_i(x, \xi_i) \) are convex and \( L_f \)-continuous respect to $x$  on \( \Omega \), where \( L_f > 0 \) is a constant.
\end{assum} \par 
\begin{assum}\label{assum3-g}
	For any $i\in \mathcal{V}$, there exists a constant ${\hat{\sigma }}>0$ s.t. $\mathbb{E}_{\xi_i}\big[ \left\| \nabla_{x} {{F}_{i}}\left(x,\xi_i\right)-\nabla f_i(x)\right\|^2 \big]\le {\hat{\sigma }^{2}}$ for all ${{x}_{i}}\in {\Omega}$.
\end{assum}\par 
\begin{rem}
Assumptions \ref{assum1-set}--\ref{assum3-g} are standard in stochastic distributed optimization \cite{gt,DSGT,noncvx-dso-zo}. Assumption \ref{assum3-g} specifically requires the stochastic gradient to have bounded second-order moments, which is reasonable in many real-world applications where data variability is limited and gradients exhibit stable statistical behavior.
\end{rem}
\subsection{Stochastic Gradient Estimation}
The intuition behind gradient estimation lies in constructing a smoothed approximation function whose unbiased gradient estimation can be constructed  using finite function evaluations. Formally, let the random vector \( u \in \mathbb{R}^d \) follow a zero-mean symmetric distribution. Define the smoothed function as $f^{\mu}(x) = \mathbb{E}_u[f\left( x + \mu u\right) ]$, where $\mu > 0$ is a parameter setting the level of smoothing.\par
In this paper, we study the two-point sampling gradient estimator, which achieves higher estimation accuracy compared to the one-point estimator. To ensure the exploration points \(x + \mu u\) remain within the constraint set \(\Omega\), we uniformly sample the perturbation direction \(u\) from \(\mathbb{S}^{d-1}\), as opposed to the Gaussian distribution used in \cite{nesterov2017, online1,b23}, which may yield unbounded values. The resulting stochastic two-point gradient estimator is defined as:
\begin{equation}\label{TP}
	\mathbf{g}^{\mu}(x) = \frac{d}{\mu} \left[\left( F(x+\mu u,\xi) - F(x,\xi) \right) u   \right]\quad \forall x \in \Omega_\epsilon,
\end{equation}
where \(\Omega_\epsilon = \{(1-\epsilon)x \mid x \in \Omega\}\) is the shrunk set with a shrinkage factor \(\epsilon > 0\), and \(d\) denotes the problem dimension. By restricting the estimator to \(\Omega_\epsilon\) and leveraging the bounded support of \(u\) on \(\mathbb{S}^{d-1}\), the query points \(x + \mu u\) are guaranteed to lie within \(\Omega\), as shown by the subsequent Lemma \ref{contracted set}.\par
\begin{lem}[\cite{flexman-onepoint}]\label{contracted set}
	For any $x\in\Omega_\epsilon$ and any unit vector $u$, $x+\mu u\in\Omega$ for any $\mu \in \left[ 0,{\epsilon} \underline{R} \right]$.
\end{lem}\par
By using the differential properties of the smoothed function, we can rigorously demonstrate that the gradient estimator \eqref{TP} provides an unbiased estimate of \(\nabla f^\mu(x)\). The following results present additional useful properties of $f^{\mu}$ and $\mathbf{g}^{\mu}$.
\begin{lem}\label{g}
	For the smoothed function $f^{\mu}$ and the gradient estimator \( \mathbf{g}^{\mu} \) as defined in \eqref{TP}, the following properties hold.\par
	1) If $f$ is convex on ${\Omega}$, then ${{f}^{\mu }}$ remains convex on $\Omega_\epsilon$.\par 
	2) If $f$ is $L_f$-Lipschitz continuous on $\Omega$, then both $f^{\mu}$ and $\nabla f^{\mu}$ are Lipschitz continuous on $\Omega_\epsilon$ with constants $L_f$ and ${dL_f}/{\mu}$, respectively. Moreover, for all $ x\in \Omega_\epsilon$, $\left| f\left(x\right)-{{f}^{\mu }\left(x\right)} \right|\le \mu {{L}_{f}}$.\par 
	3)	If $f$ is $L_f$-Lipschitz continuous on $\Omega$, $\mathbb{E}_{u\in\mathbb{S}^{d}}\left[ \left\| \mathbf{g}^\mu\left( x \right) \right\|\right]\le dL_f$ for all $ x\in \Omega_\epsilon$.\par
	4) If \( f \) has \( L_m \)-Lipschitz continuous gradient on \( \Omega \), then, for all $ x\in \Omega_\epsilon$, $\left\| \nabla {{f}^{\mu}_{i}}\left(x\right) -\nabla f_i\left( x \right) \right\|\leq \mu L_m$ and $\mathbb{E}_{u\in\mathbb{S}^{d}}\big[ \left\| \mathbf{g}^\mu\left( x \right) \right\|^2\big]\leq d\left\| \nabla F\left(x,\xi \right)  \right\|^2+d^2\mu^2L_m^2/2$.
\end{lem}\par
\begin{IEEEproof}
	For the proofs of Lemma \ref{g}, please refer to \cite{online2}, Lemma 5 in \cite{zero-ordernonconvex} and Proposition 7.6 in \cite{admm1}.
\end{IEEEproof}

\subsection{Compression operators}
This subsection introduces the conditions on the compression operator $\mathcal{C}\left(\cdot \right) :\mathbb{R}^d\to\mathbb{R}^d$.\par 
\begin{assum}\label{assum-c}
	The rescaled compressor \( \mathcal{C}_r\left(\cdot \right): \mathbb{R}^d \to \mathbb{R}^d \), defined as \( \mathcal{C}_r(x) = \frac{1}{r}\mathcal{C}(x) \) with a scaling factor \( r \geq 1 \), satisfies
		\begin{equation}
		\mathbb{E}_{\mathcal{C}} \left[  \left\| {\mathcal{C}_r\left(x \right) } - x \right\|^2 \right]  \leq (1 - \omega) \| x \|^2
	\end{equation}
	for all \( x \in \Omega \), where \( \omega \in (0,1] \) is a positive constant setting the compression level and $\mathbb{E}\left[ \cdot \right] $ is taken over the internal randomness of $\mathcal{C}\left(\cdot \right)$.
\end{assum}\par 
The bounded relative error condition for compression operators, standard in compressed communication frameworks \cite{b25,b26,unbiased1,Choco-gossip,onoine-choco}, ensures controlled compression error. Departing from classical contraction assumptions, our rescaled stochastic compressors generalize to some non-contractive operators, including norm-sign compressors \cite{b33}.
While unbiased compression is commonly used in DCO \cite{b25,b26,unbiased1}, this assumption precludes deterministic schemes like Top-$k$ sparsification, motivating our focus on biased yet practical compressors.\par 
To quantify communication efficiency, consider transmitting a $d$-dimensional real vector $x$ (64-byte representation).\par 
\begin{exmp}[Norm-Sign Compressor]
	The operator \( \mathcal{C}_1\left(\cdot \right) : \mathbb{R}^d \to \mathbb{R}^d \) defined by
\begin{equation*}
	\mathcal{C}_{1}\left( x\right) =\frac{\left\|x\right\|_\infty}{2}\text{sgn}\left( x\right) \text{,}
\end{equation*}
	where \( \|x\|_\infty := \max_{1\leq i\leq d} |x_i| \) and $\text{sgn}$ is the signum function. Note that $	\mathcal{C}_{1}\left( \cdot\right)$ is a biased and non-contractive compressor, but it satisfies Assumption \ref{assum-c} with $r={1}/{d}$ and $\omega={1}/{d^{2}}$. Transmitting $	\mathcal{C}_1(x)$ requires \( d \)-bit (signs) + 64-bit (scalar \( \|x\|_\infty/2 \)), which totals \( d + 64 \) bits. Thus, this compressor trades bias for 32$\times$ bandwidth reduction compared to 32-bit float transmission.
\end{exmp}
\section{Distributed Stochastic Zeroth-Order Algorithms with Compressed Communication}
To reduce communication overhead, a quantized version of DGD, Choco-Gossip algorithm \cite{Choco-gossip} firstly achieved sublinear convergence for a class of biased but $\delta$-contracted compressors. The core idea of Choco-Gossip lies in the introduction of a reference variable $a_{i,k}$. The difference between the state vector $x_{i,k}$ and the reference vector $a_{i,k}$, termed the “innovation vector”, is compressed and added back to $x_{i,k}$ to yield the updated estimate $\hat{x}_{i,k+1}$. The algorithm proceeds as follows: 
\begin{subequations}
	\begin{align}
		&y_{i,k+1}=x_{i,k}-\eta_k\nabla f_{i}\left(x_{i,k} \right) \label{bb1}\\
		&x_{i,k+1}=y_{i,k+1}+\gamma\left(\sum\nolimits_{j\in\mathcal{N}_i}\mathbf{W}_{i,j}\hat{x}_{j,k}-\hat{x}_{i,k} \right)\label{bb2}\\
		&\hat{x}_{i,k+1}=a_{i,k}+\mathcal{C}\left(x_{i,k}-a_{i,k} \right)\label{bb3}
	\end{align}
\end{subequations}
where $\eta_k$ and $\gamma$ are positive parameters.\par
Although the algorithm reduces the communication burden, the agents still need to transmit the vector value \( a_{i,k} \) at the time of execution due to the summation term. In recent work \cite{onoine-choco}, \cite{b33} and \cite{nonconvexcom}, only the compression vector \( q_{i,k} \) needs to be exchanged between agents. However, each agent \( i \) must additionally maintain and update \( \mathcal{N}_i \) auxiliary variables, which is highly resource-intensive. Inspired by this, we propose the communication-efficient distributed zeroth-order (Com-DSZO) algorithm for more general compression operators, as summarized in Algorithm 1.\par 
\begin{algorithm}
	\caption{ Com-DSZO: Communication-efficient Distributed Stochastic Zeroth-Order algorithm}
	\par	\textbf{Input:} $\eta_k$, $\mu$, ${\epsilon}$, $\psi$, and $\gamma$.
	\par	\textbf{Initialize:} $x_{i,0}\in\Omega, \hat{x}_{i,0} = {b}_{i,0}= \mathbf{0}_d$, and $q_{i,0}=\mathcal{C}\left(x_{i,0}\right)$, $\forall i\in\mathcal{V}$.
	\par	\textbf{for} {$k = 0,1,\ldots$ \textbf{do}} 
	\par \quad	\textbf{for} {$i = 1$ \textbf{to} $n$ in parallel \textbf{do}} 
	\par \qquad	Generate $u_{i,k}$ , $\xi_{i,k}$ and construct $\mathbf{g}_i^\mu(x_{i,k})$ by \eqref{TP}.
	\par \qquad
	Send $q_{i,k}$ to $j\in\mathcal{N}_i$, receive $q_{j,k}$ form $j\in\mathcal{N}_i$.\par \qquad
	Update
	\begin{align}
		&\hat{x}_{i,k+1} =\hat{x}_{i,k} + \psi q_{i,k}\label{a1}\\
		&b_{i,k+1}=b_{i,k}+\psi\left(q_{i,k}-\sum\nolimits _{j\in\mathcal{N}_i}\mathbf{W}_{i,j}q_{j,k}\right)\label{a2} \\
		&x_{i,k+1}= \Pi_{\Omega_\epsilon}\left[  x_{i,k} -\gamma b_{i,k+1} - \eta_k \mathbf{g}_{i}^{\mu}\left(x_{i,k}\right)\right]\label{a3} \\
		&q_{i,k+1} = \mathcal{C}(x_{i,k+1} - \hat{x}_{i,k+1})\label{a4}
	\end{align}
	\par \quad \textbf{end for}
	\par \textbf{end for}	
\end{algorithm}
In this algorithm, each agent $i$ maintains and updates four variables: one state variable $x_{i,k}$ and three auxiliary variables $\hat{x}_{i,k}$, $b_{i,k}$, and $q_{i,k}$. The variable  $\hat{x}_{i,k}$ acts as a reference point for $x_{i,k}$. The difference between $x_{i,k}$ and the auxiliary variable $\hat{x}_{i,k}$ is compressed in \eqref{a4}, and then added back to $x_{i,k}$ in \eqref{a1} to obtain $x_{i,k+1}$. As \(\hat{x}_{i,k}\) gradually approaches \(x_{i,k}\), the compression error decreases to $0$, as stated in Assumption 5.\par
The auxiliary variable $b_{i,k}$ can be considered as a local backup of neighboring information, eliminating the need for storing additional estimates of the state of the neighbors. Notice that $\hat{x}_{i,0} = {b}_{i,0}= \mathbf{0}_d$ in the initial step,   from \eqref{a1} and \eqref{a2}, we can derive that ${b}_{i,k}=\hat{x}_{i,k}-\sum\nolimits_{j\in\mathcal{N}_i}\mathbf{W}_{i,j}\hat{x}_{j,k}$ by simple mathematical induction, which implies that \eqref{a3} can be rewritten in a form similar to \eqref{bb2} in Choco-Gossip algorithm. The low-bit variable $q_{i,k}$ is exchanged between neighboring agents. The agent performs projected gradient descent on the shrunk set to obtain the updated state as shown in \eqref{a3}.\par 
\section{Main Results}
This section focuses on establishing the main convergence results for the Com-DSZO algorithm. Before proceeding, the following notation is introduced.\par 
Define ${{\tilde{x}}_{i,k+1}}=x_{i,k} -\gamma b_{i,k+1} - \eta_k \mathbf{g}_{i}^{\mu}\left(x_{i,k}\right)$, the projection error ${{\varphi }_{i,k}}={{x}_{i,k}}-{{\tilde{x}}_{i,k}}$, and the network-averaged state $\bar{x}_k\triangleq\frac{1}{n}\sum\nolimits_{i=1}^nx_{i,k}$.
	Denote $\mathbf{X} := [x_1, x_2, \dots, x_n]^\mathsf{T} \in \mathbb{R}^{n \times d}$, $\mathbf{\overline{X}} :=  \frac{1}{n}\mathbf{1}\mathbf{1}^\mathsf{T}\mathbf{X} \in \mathbb{R}^{n \times d}$, $\mathbf{G}^{\mu} := [\mathbf{g}_{1}^{\mu}, \mathbf{g}_{2}^{\mu}, \dots, \mathbf{g}_{n}^{\mu}]^\mathsf{T}\in\mathbb{R}^{n\times d}$ and $\mathbf{\Phi} := [\varphi_1, \varphi_2, \dots, \varphi_n]^\mathsf{T} \in \mathbb{R}^{n \times d}$
	. At the \( k \)-th iteration, their values are denoted by \( \mathbf{X}_k, \mathbf{\overline{X}}_k \), $\mathbf{G}_k^{\mu}$ and $\mathbf{\Phi}_k$, respectively. Other auxiliary variables of the agents (in compact form), \( \mathbf{\widehat{X}}_k, \mathbf{B}_k, \mathbf{Q}_k, \mathbf{\widetilde{X}}_k \) are defined similarly. To facilitate the convergence analysis, we denote $\mathbf{g}_{i,k}^{\mu} \triangleq \mathbf{g}_i^\mu(x_{i,k})$ as the local gradient estimator at iteration $k$. Define two error metrics: consensus error $U_{1,k} \triangleq \|\mathbf{X}_k - \overline{\mathbf{X}}_k\|_F^2$ and compression error $U_{2,k} \triangleq \|\widehat{\mathbf{X}}_{k+1} - \mathbf{X}_k\|_F^2$, with composite dynamics $\mathbf{e}_k \triangleq \mathbb{E}[U_{1,k} + U_{2,k}]$.\par 
	For clarity, define the compression operator across agents as $\mathcal{{C}}\left({\mathbf{X}_{k}}\right):= [\mathcal{C}\left(x_{1,k}\right),\ldots,\mathcal{C}\left(x_{n,k}\right)]^\mathsf{T}$. Let $\Omega_{\epsilon}^n$ represents the Cartesian product $\Omega_{\epsilon} \times \ldots \times \Omega_{\epsilon}$ across $n$ components. Algorithm 1 can then be compactly expressed as follows:
	\begin{subequations}
		\begin{alignat}{4}
			&{{\mathbf{\widehat{X}}}_{k+1}}={{\mathbf{\widehat{X}}}_{k}}+\psi {\mathbf{Q}_{k}}\label{a} \\
			&\mathbf{B}_{k+1}=\mathbf{B}_k+\psi\left(\mathbf{I}-\mathbf{W} \right)\mathbf{Q}_k\label{b}\\  
			& {\mathbf{X}_{k+1}}=\Pi_{\Omega_{\epsilon}^n}{\left[{\mathbf{X}_{k}}-\gamma\mathbf{B}_{k+1}-{{\eta }_{k}}\mathbf{G}_{k}^{\mu}\right]} \label{c} \\
			& {\mathbf{Q}_{k+1}}=\mathcal{C}\big( \mathbf{X}_{k+1}-{{\mathbf{\widehat{X}}}_{k+1}} \big)\label{d} 
		\end{alignat}
\end{subequations}
\subsection{Convergence analysis for smooth functions}
This subsection develop the convergence analysis of the Com-DSZO algorithm for stochastic constraint DCO with smooth functions. That is, under the following assumption.\par
\begin{assum}\label{smooth}
	Each ${{f}_{i}}$ has $L_{{m}}$-Lipschitz gradient, i.e., $\left\| \nabla {{f}_{i}}\left( x \right)-\nabla {{f}_{i}}\left( y \right) \right\|\le L_{{m}}\left\| x-y \right\|$ for any $x,y\in \Omega $.
\end{assum}\par
\begin{thm}\label{result 1 for smooth function}
	Consider the stochastic constrained DCO in \eqref{1}. Under Assumptions \ref{assum G}-\ref{smooth}, let $\{x_{i,k}\}_{k>0}$ be generated by Com-DSZO with parameters:  
	$\mu \!=\!{1}/{\sqrt{T}}$, $\epsilon\!=\!{\mu}/{2\underline{R}}$, $\psi \in (0, r^{-1}]$, $\eta_k\!=\!{1}/{(M_1\sqrt{k+c})}$, where $M_1 = \left(5+{10}/{n}\right)4dL_m$, $c>4$ and $\gamma$ is defined in \eqref{gamma}. Then for any optimal point $x^* \in \Omega$:
	\begin{align}\label{thm1}
		&\sum\limits_{k=1}^T \mathbb{E}\big[f(\bar{x}_k) - f(x^*)\big] 
		\leq \mathcal{O}(1) + \mathcal{O}(\sqrt{T}) + \mathcal{O}\left( \frac{1}{T}\right)  \nonumber\\
		&+ \mathcal{O}\left( \frac{1}{\sqrt{T}}\right)  + \mathcal{O}\left( \frac{\kappa_3}{1-\kappa_0} \left( \sqrt{T+c} + \ln(T+c)\right)  \right)\text{,}
	\end{align}
where $\kappa_0\leq\max\left(1,\mathbf{C}_1,\mathbf{C}_2\right)$ and $\kappa_3=nd^2L_f^2\mathbf{\widetilde{C}}_3$. The constants $\mathbf{C}_1$, $\mathbf{C}_2$, $\mathbf{\widetilde{C}}_3$ (dependent on dimension $d$, compression resolution $\omega$ and spectral radius $\delta$) are formally defined in Appendix A. 
\end{thm}\par 
Theorem \ref{result 1 for smooth function} is proved via the following lemma establishing bounds for disagreement and compression error accumulation. The proof of Lemma \ref{smooth-e} is provided in Appendix A.
\begin{lem}\label{smooth-e}
	Under the same conditions as in Theorem \ref{result 1 for smooth function}, we have
\begin{align}\label{lem-e}
	\sum\limits_{k=1}^T\mathbf{e}_{k}&\leq \frac{\kappa_3}{1-\kappa_0}\sum\limits_{k=0}^{T-1}\eta_{k}^2\text{.}
\end{align}
	Moreover, let $0<\gamma\leq \frac{\delta c_0^3}{(\delta+2)c_0^2\beta^2+(c_0+1)^2\beta^2}$, and $\kappa_0=1-c_0$ with $c_0=\frac{\psi r\omega}{2}$. Then, we have
\begin{align}\label{e-trade-off}
	\sum\limits_{k=1}^T\mathbf{e}_{k}&\leq\mathcal{O}\left(\frac{d^2}{\left( \psi r \omega\right) ^2 }\sum\nolimits_{k=0}^{T-1}\eta_{k}^2\right) \text{.}
\end{align}\par 
\end{lem} \par
\begin{IEEEproof}[Proof of Theorem 1]
	We start by examining the evolution of the optimization error $V_{k+1}\left({\epsilon} \right)=\left\|{\bar{x}}_{k+1}-\left(1-{\epsilon} \right)x^*\right\|^2$ for any $x^*\in\Omega$. Recalling for the definition of $\varphi_k$ and \eqref{a3}, it follows that
	\begin{align}\label{thm1-1}
		V_{k+1}\left({\epsilon}\right)&= \left\|{\bar{\tilde{x}}}_{k+1}-\left(1-{\epsilon} \right)x^*+\bar{\varphi}_{k+1}\right\|^2\nonumber\\
		&=\left\|{\bar{x}}_{k}-\left(1-{\epsilon} \right)x^*-\eta_k\mathbf{\bar{g}}_{k}^{\mu}\right\|^2+\left\|\bar{\varphi}_{k+1}\right\|^2\nonumber\\
		&\quad+2\left\langle{\bar{x}}_{k}-\left(1-{\epsilon} \right)x^*-\eta_k\mathbf{\bar{g}}_{k}^{\mu},\bar{{\varphi}}_{k+1}\right\rangle\text{,} 
	\end{align}
where $\bar{\varphi}_{k+1}=\frac{1}{n}\sum\nolimits_{i=1}^n\varphi_{i,k+1}$, the second equality holds since the network averaging residual $\bar{b}_{k+1}=\frac{1}{n}\sum\nolimits_{i=1}^n[\hat{x}_{i,k+1}-\sum\nolimits_{j\in\mathcal{N}_i}\mathbf{W}_{i,j}\hat{x}_{j,k+1}]=0$ under the doubly stochasticity of $\mathbf{W}$. We now proceed to provide the bounds for each term on the right-hand side of \eqref{thm1-1} in expectation. For the first term,
\begin{align}\label{thm1-2}
	&\left\|{\bar{x}}_{k}-\left(1-{\epsilon} \right)x^*-\eta_k\mathbf{\bar{g}}_{k}^{\mu}\right\|^2 \nonumber\\
	&=V_k(\epsilon)+\eta_k^2\big\|\mathbf{\bar{g}}_{k}^{\mu}\big\|^2-2\eta_k\left\langle{\bar{x}}_{k}-\left(1-{\epsilon} \right)x^*,\mathbf{\bar{g}}_{k}^{\mu} \right\rangle .
\end{align}
By the definition of $\mathbf{\bar{g}}_k^{\mu}=\frac{1}{n}\sum\nolimits_{i=1}^n\mathbf{g}_{i,k}^{\mu}$, it follows that 
$ \big\|\mathbf{\bar{g}}_{k}^{\mu}\big\|^2\leq \frac{1}{n} \left\|\mathbf{G}_{k}^{\mu}\right\|^2$.
Moreover, Assumption \ref{smooth} gives that
	\begin{align}
	&\mathbb{E}\left[ \left\langle {\bar{x}}_{k}-\left(1-{\epsilon} \right) x^*,\mathbf{\bar{g}}_k^{\mu}\right\rangle\right]\nonumber\\
	&=\mathbb{E}\left[ \left\langle {\bar{x}}_{k}-x_{i,k},\overline{\nabla}{f}^{\mu}_k\right\rangle\right]+\mathbb{E}\left[\left\langle x_{i,k}-\left(1-{\epsilon} \right) {x}^*,\overline{\nabla}{f}^{\mu}_k \right\rangle\right] \nonumber\\
	&\geq \mathbb{E}\left[ f\left(\bar{x}_{k} \right) -f \left(\left(1-{\epsilon} \right){x}^* \right)\right]-2L_f\mu-\frac{L_{m}}{n} \mathbb{E}\left[U_{1,k}\right]\nonumber\\
	&\geq \mathbb{E}\left[ f\left(\bar{x}_{k} \right) -f \left({x}^* \right)\right]-L_f\bar{R}\epsilon-2L_f\mu-\frac{L_{m}}{n} \mathbb{E}\left[ U_{1,k}\right]\text{,}\nonumber
\end{align}
where $\overline{\nabla}{f}_{k}^{\mu}\!=\!\frac{1}{n}\sum\nolimits_{i=1}^n\nabla f_{i}^\mu$, the equality holds by $\mathbb{E}\left[\mathbf{g}_i^{\mu} \right]\!=\!\nabla f_i^{\mu}$; the fist inequality is derived from the convexity and $L_m$-smoothness of $f_i$ and Lemma \ref{g}-2); and the final inequality applies the $L_f$-Lipschitz continuity of $f_{i}$ and Assumption \ref{assum1-set}. Substituting the above two inequalities into \eqref{thm1-2} yields that
\begin{align}\label{t1}
	&\mathbb{E}\left[ \left\|{\bar{x}}_{k}-\left(1-{\epsilon} \right)x^*-\eta_k\mathbf{\bar{g}}_{k}^{\mu}\right\|^2\right]\nonumber\\
	&\leq \mathbb{E}\left[ V_k(\epsilon)\right] -2\eta_k\mathbb{E}\left[ f\left(\bar{x}_{k} \right) -f \left({x}^* \right)\right]+\frac{\eta_k^2}{n}\mathbb{E}\left[ \left\|\mathbf{G}_k^{\mu} \right\|_F^2 \right] \nonumber\\
	&\quad+\left( \bar{R}\epsilon+2\mu\right)2L_f\eta_k+\frac{2L_m\eta_k}{n}\mathbb{E}\left[U_{1,k}\right]. 
\end{align}\par 
The non-expansiveness of the projection operator $\Pi_{\Omega_\epsilon}$ establishes $\mathbb{E}\left\|\Phi_{k+1}\right\|_F^2\leq \mathbb{E}\big\|\mathbf{\tilde{X}}_{k+1}-\mathbf{X}_k\big\|_F^2$. Setting $\alpha_4=1$ in \eqref{3-2} then gives
\begin{align}\label{t2}
	{\mathbb{E}}\left\| \bar{\varphi}_{k+1}\right\|^2  	&\leq\frac{1}{n}\mathbb{E}\left\|\Phi_{k+1}\right\|_F^2\nonumber\\
	&\leq \frac{2\gamma^2\beta^2}{n}\left(  2\mathbb{E}\left[U_{1,k}\right]+\mathbb{E}\left[U_{2,k}\right]\right) +\frac{4\eta_k^2}{n}\mathbb{E}\left[\left\|\mathbf{G}_k^{\mu}\right\|_F^2\right]  \nonumber\\
	&\leq \frac{4\gamma^2\beta^2}{n}\mathbf{e}_k+\frac{4\eta_k^2}{n}\mathbb{E}\left[\left\| \mathbf{G}_{k}^{\mu}\right\| _F^2\right] .
\end{align}\par 
For the last term of \eqref{thm1-1}, given that $\left\langle y-x, \Pi_{\Omega_\epsilon}\left[y\right] -y\right\rangle \leq 0$ for any $x\in \Omega_\epsilon$ and $y\in \mathbb{R}^d$, we have
\begin{align}
	\left\langle{\bar{x}}_{k}-\left(1-{\epsilon} \right)x^*,\bar{{\varphi}}_{k+1}\right\rangle
	&\leq \frac{1}{n}\sum\nolimits_{i=1}^{n}\left\langle{\bar{x}}_{k}-\tilde{x}_{i,k+1},\bar{{\varphi}}_{k+1} \right\rangle\nonumber\\
	&=-\eta_k \left\langle\mathbf{\bar{g}}_k^{\mu},\bar{{\varphi}}_{k+1} \right\rangle\text{,}\nonumber
\end{align}
which further implies that
\begin{align}\label{t3}
	&2\left\langle{\bar{x}}_{k}-\left(1-{\epsilon} \right)x^*-\eta_k\mathbf{\bar{g}}_{k}^{\mu},\bar{{\varphi}}_{k+1}\right\rangle\leq -4\eta_k \left\langle\mathbf{\bar{g}}_k^{\mu},\bar{{\varphi}}_{k+1} \right\rangle\nonumber\\
	&\leq \frac{4\eta_k}{n^2}\left\|\mathbf{G}_k^{\mu} \right\|\cdot\left\|\mathbf{\Phi}_{k+1} \right\|\leq\frac{2\eta_k^2}{n^2}\left\|\mathbf{G}_k^{\mu} \right\|^2+\frac{2}{n^2}\left\|\mathbf{\Phi}_{k+1} \right\|^2\text{,}
\end{align}
where the final inequality holds by the Cauchy-Schwarz inequality. 
By combining \eqref{thm1-1}--\eqref{t3} yields 
\begin{align}\label{thm1-4}
	&\mathbb{E}\left[V_{k+1}\left({\epsilon}\right)\right] \leq \mathbb{E}\left[ V_k(\epsilon)\right] -2\eta_k\mathbb{E}\left[ f\left(\bar{x}_{k} \right) -f \left({x}^* \right)\right] \nonumber\\
	&+\left( \frac{10}{n^2}+\frac{5}{n}\right) {\eta_k^2}\mathbb{E}\left[\left\| \mathbf{G}_{k}^{\mu}\right\| _F^2\right]+\left( \frac{2}{n^2}+\frac{1}{n}\right) {4\gamma^2\beta^2}\mathbf{e}_k\nonumber\\
	&+\left( \bar{R}\epsilon+2\mu\right)2L_f\eta_k+\frac{2L_m\eta_k}{n}\mathbb{E}\left[U_{1,k}\right].
\end{align}
Given that Lemma \ref{g} and Assumption \ref{assum3-g}, we have
\begin{align}\label{smooth-G}
    &\mathbb{E}\left[ \left\|\mathbf{G}_{k}^{\mu}\right\|_F^2\right]= \sum\nolimits_{i=1}^n\mathbb{E}\Big[ \big\|\mathbf{g}_{i,k}^{\mu}\big\|^2\Big]\nonumber\\
	&\leq \hspace{-2pt}2nd{{\hat{\sigma} }^{2}}\hspace{-2pt}+\hspace{-2pt}2d\sum\nolimits_{i=1}^n\mathbb{E}\left[ \left\| \nabla f_i\left( x_{i,k}\right) \right\|^2\right] \hspace{-2pt}+\hspace{-2pt}\frac{nd^2L_{m}^{2}\mu^2}{2} \text{.} 
\end{align}
Using the smoothness of the function $f_i$ yields
\begin{align}
	& {{{\left\| \nabla {{f}_{i}}\left( {{x}_{i,k}} \right) \right\|}^{2}}}
	\hspace{-2pt}\le\hspace{-2pt}2{{L_{m}^{2}}{{\left\| {{x}_{i,k}}-{{{\bar{x}}}_{k}} \right\|}^{2}}}\hspace{-2pt}+\hspace{-2pt}4{L_{m}}\left(f_i\left( {{{\bar{x}}}_{k}} \right)\hspace{-2pt}-\hspace{-2pt}{{f}_{i}}\left( {{x}^{*}} \right)  \right).\nonumber
\end{align}
By combining the above two inequalities, we can infer that
\begin{align}\label{thm1-3}
	\mathbb{E}\left[ \big\|\mathbf{G}_{k}^{\mu}\big\|^2\right]&\leq 2nd{{\hat{\sigma} }^{2}}+ \frac{nd^2{{L_{m}^{2}\mu }^{2}}}{2}+{4dL_{m}^{2}}\mathbb{E}\left[U_{1,k}\right]  \nonumber\\
	&\quad+ {8nd{L_{m}}}\left(f\left( {{{\bar{x}}}_{k}} \right)- {f}\left( {{x}^{*}} \right)  \right).
\end{align}
Substituting \eqref{thm1-3} into \eqref{thm1-4} gives that
\begin{align}\label{thm1-5}
	&\mathbb{E}\left[V_{k+1}\left({\epsilon}\right)\right] \leq \mathbb{E}\left[ V_k(\epsilon)\right] -2\eta_k\left[1-M_1\eta_k\right] \mathbb{E}\left[ f\left(\bar{x}_{k} \right) -f \left({x}^* \right)\right] \nonumber\\
	&+\left( \frac{2}{n^2}\hspace{-2pt}+\hspace{-2pt}\frac{1}{n}\right) {4\gamma^2\beta^2}\mathbf{e}_k\hspace{-2pt}+\hspace{-2pt}\left( \frac{2L_m\eta_k}{n}\hspace{-2pt}+\hspace{-2pt}4dL_{m}^{2}M_0 {\eta_k^2}\right)\mathbb{E}\left[U_{1,k}\right]\nonumber\\
	&+\left( \bar{R}\epsilon+2\mu\right)2L_f\eta_k+\left(4nd{{\hat{\sigma} }^{2}}+ {nd^2L_m^2\mu^{2}}\right)\frac{M_0{\eta_k^2}}{2}\text{,}
\end{align}
where $M_0=\left( \frac{10}{n^2}+\frac{5}{n}\right)$ and $M_1=4ndL_mM_0$. Note that
\begin{align}\label{1-2}
	&\sum\nolimits_{k=1}^T\frac{1}{z_k}\left\lbrace \mathbb{E}\left[  V_{k}\left({\epsilon} \right)\right] -\mathbb{E}\left[ V_{k+1}\left({\epsilon} \right)\right]\right\rbrace \nonumber\\
	&=\frac{\mathbb{E}\left[  V_{1}\left({\epsilon} \right)\right] }{z_1}\hspace{-2pt}-\hspace{-2pt}\frac{\mathbb{E}\left[  V_{T+1}\left({\epsilon} \right)\right]}{z_T}\hspace{-2pt}+\hspace{-2pt}\sum\nolimits_{k=2}^T\hspace{-2pt}\left(\hspace{-1pt}\frac{1 }{z_k}\hspace{-2pt}-\hspace{-2pt}\frac{1}{z_{k-1}}\hspace{-1pt}\right)\hspace{-2pt}{\mathbb{E}\left[  V_{k}\left({\epsilon} \right)\right]}\nonumber\\
	&\leq \frac{4\bar{R}^2}{z_1}+4\bar{R}^2\left(\frac{1 }{z_T}-\frac{1}{z_{1}} \right)=\frac{4\bar{R}^2}{z_{T}} \text{,}
\end{align}
where the fact $\left\|x-y\right\|^2\leq 4\overline{R}^2 $ for any $x,y\in \Omega$ is used. Let $\eta_k=\frac{1}{M_1(k+c)^a}$ with $c>0$. By rearranging the terms of \eqref{thm1-5} and summing both sides over $k$ from $1$ to $T$, it follows that
\begin{align}\label{thm1-7}
	&\sum\nolimits_{k=1}^T\mathbb{E}\left[ f\left(\bar{x}_{k} \right) -f \left({x}^* \right)\right]\leq \frac{2\bar{R}^2M_1(T+c)^{2a}}{(T+c)^a-1}\nonumber\\
	&+\left( \frac{2}{n^2}+\frac{1}{n}\right) {2\gamma^2\beta^2M_1}\cdot\sum\nolimits_{k=1}^T\frac{(k+c)^{2a}}{(k+c)^a-1}\mathbf{e}_k\nonumber\\
	&+\left( \bar{R}\epsilon+2\mu\right)L_f\cdot\sum\nolimits_{k=1}^T\frac{(k+c)^{a}}{(k+c)^a-1}+ \frac{L_m}{n}\cdot\sum\nolimits_{k=1}^T\mathbf{e}_k\nonumber\\
	&+\left( \frac{L_m}{n}+2dL_{m}^{2}\frac{M_0}{M_1}\right) \cdot\sum\nolimits_{k=1}^T\frac{\mathbf{e}_k}{(k+c)^a-1}\nonumber\\
	&+\frac{4{{\hat{\sigma} }^{2}}+ {dL_m^2\mu^{2}}}{16L_m}\cdot\sum\nolimits_{k=1}^T\frac{1}{(k+c)^a-1}.
\end{align}\par 
Given \( c^a > 2 \), the inequalities \( \frac{(k+c)^{2a}}{(k+c)^a-1}\mathbf{e}_k \leq \big(2+(k+c)^a\big)\mathbf{e}_k \) holds. Lemma \eqref{lem-e} further provides the summations \( \sum_{k=1}^T\mathbf{e}_k \leq \frac{\kappa_3}{1-\kappa_0}\sum\nolimits_{k=1}^T\eta_k^2 \). From \( c^a > 2 \), we deduce \( \frac{1}{(k+c)^a -1} < \frac{2}{(k+c)^a} \) for all \( k \geq 0 \), leading to \( \sum_{k=0}^T \frac{1}{(k+c)^a -1} < 2(T+c)^{1-a} \) and $\sum\nolimits_{k=1}^T\frac{\mathbf{e}_k}{(k+c)^a-1}< \sum\nolimits_{k=1}^T\mathbf{e}_k$. On the other hand, 
from the recursive relationship in \eqref{e-diedai}, we derive the error propagation dynamics $\frac{\mathbf{e}_{k+1}}{\eta_{k+1}} \leq \kappa_0\frac{\mathbf{e}_k}{\eta_k}\cdot\frac{\eta_k}{\eta_{k+1}} + \mathbf{C}_3\eta_k\mathbb{E}\left[\|\mathbf{G}_k^{\mu}\|_F^2\right]\cdot\frac{\eta_k}{\eta_{k+1}}$.
Following analogous arguments to \eqref{e2} with $\mathbf{e}_0 = 0$, the error sequence admits: $\frac{\mathbf{e}_k}{\eta_k} \leq \kappa_3\cdot\frac{\eta_{k-1}^2}{\eta_k^2}\sum_{\tau=0}^k\kappa_0^{k-1-\tau}\eta_k$.
For the step size $\eta_k = \frac{1}{M_1\sqrt{k+c}}$, observe that $\frac{\eta_{k-1}^2}{\eta_k^2} = \frac{k+c}{k-1+c} < 2$ for any $k\geq 1$, where the inequality holds since $c>1$. Summation over $k=1$ to $T$ gives
\begin{align}\label{thm1-6}
	\sum\nolimits_{k=1}^T\frac{\mathbf{e}_k}{\eta_k} &\leq \frac{2\kappa_3}{1-\kappa_0}\sum\nolimits_{k=0}^{T-1}\eta_k.
\end{align}
Substituting these into \eqref{thm1-7} with parameters \( \mu = \frac{1}{\sqrt{T}} \) and \( \epsilon = \frac{2\underline{R}}{\sqrt{T}} \), we can derive:
\begin{align}\label{thm1-9}
	&\sum\nolimits_{k=1}^T\mathbb{E}\left[ f\left(\bar{x}_{k} \right) -f \left({x}^* \right)\right]\leq {2\overline{R}^2M_1}\left(2+{\sqrt{T+c}}\right) \nonumber\\
	&+ \left( \frac{2}{n^2}+\frac{1}{n}\right) \frac{4\gamma^2\beta^2\kappa_3}{M_1(1-\kappa_0)}\cdot \ln(T+c) \nonumber\\
	&+\left( \frac{2L_m}{n}+2dL_{m}^{2}\frac{M_0}{M_1}\right) \frac{\kappa_3}{M_1^2(1-\kappa_0)} \cdot \ln(T+c) \nonumber\\
	&+\left[\left(\frac{2}{n^2}+\frac{1}{n}\right) \frac{2\gamma^2\beta^2\kappa_3}{M_1(1-\kappa_0)}+\frac{\hat{\sigma}^2}{2L_m}\right] \cdot\sqrt{T+c} \nonumber\\
	&+\left( 2+{\frac{\overline{R}}{2\underline{R}}}\right)L_f\cdot\left(\sqrt{T}+2+\frac{2\sqrt{c}}{\sqrt{T}}\right)\nonumber\\
	&+\frac{dL_m}{8}\left(\frac{1}{\sqrt{T}}+\frac{\sqrt{c}}{T} \right) .
\end{align} 
We conclude the proof of Theorem 1.
\end{IEEEproof}
\begin{corollary}
	Under the same conditions as in Theorem \ref{result 1 for smooth function}. Let $\kappa_0=1-c_0$ with $c_0=\frac{\psi r\omega}{2}$. Then for any optimal point $x^* \in \Omega$:
	\begin{align}\label{corollary1}
		&\sum\limits_{k=1}^T\mathbb{E}\left[ f\left(\bar{x}_{k} \right) -f \left({x}^* \right)\right]\leq  \mathcal{O}\left(\frac{d^2\sqrt{T+c}}{\left( \psi r\omega\right) ^2}\right). 
	\end{align}\par 
	Under high compression fidelity (\(\omega \to 1\), \(\psi r \to 1\)) yielding \(\psi r\omega \approx 1\), we derive  
	$\sum_{k=1}^T \mathbb{E}\left[ f\left( \bar{x}_k \right) - f\left( x^* \right) \right] < \mathcal{O}\left( \sqrt{T+c} \right)$, which matches the convergence rate of stochastic distributed ZO method using exact first-order information \cite{sundharram_distributed_2010}.
\end{corollary}\par
\begin{IEEEproof}
The proof technique of Lemma~3 establishes that the condition $\kappa_0 < \max(\mathbf{C}_1,\mathbf{C}_2)$ holds when $\kappa_0 = 1 - c_0$. This parameter selection yields $1 - \kappa_0 = \frac{\psi r\omega}{2}$. Furthermore, observe that $\kappa_3 = nd^2L_f^2\widetilde{\mathbf{C}}_3 = \mathcal{O}(d^2c_0^{-1})$. Substituting these quantified relationships into the convergence result in Theorem 1 yields the desired convergence-compression trade-off.
\end{IEEEproof}
\begin{rem}
First, the Com-DSZO algorithm achieves the well-established convergence rate of $\mathcal{O}(1/\sqrt{T})$ for solving \eqref{1} (see Theorem \ref{thm1}), matching the convergence rates achieved by distributed algorithms employing exact communication (e.g., [12], [15], [19]–[20]). This demonstrates that our method preserves convergence guarantees while significantly improving communication efficiency in stochastic distributed constrained optimization. Second, Theorem 1 provides the first convergence proof for gradient-free optimization in the presence of both compression and stochasticity. Compared to existing compressed first-order methods [25]–[30], [32]–[35], our gradient-free approach introduces additional perturbations due to the inherent bias in the gradient estimator, necessitating a more intricate analysis to address these compounded challenges in distributed stochastic settings. Finally, our analysis extends to general compressors beyond the commonly studied unbiased and c-contracting cases [27]–[28], [30]–[31], [33], [36]. Notably, Theorem 1 establishes the first explicit characterization of the compression-convergence trade-off for compressed stochastic optimization under this broader compressor class.
\end{rem}

\subsection{Convergence analysis for non-smooth functions}
In this section, we establish the convergence properties of the proposed Com-DSZO algorithm under the scenario where the local objective functions $f_i$ in problem \eqref{1} are convex but not necessarily smooth.\par 
\begin{thm}\label{result 1 for non-smooth function}
	Consider optimization problem (1). Suppose that Assumptions \ref{assum G}--5 hold and ${{\left\{ {{x}_{i,k}} \right\}}_{k> 0}}$ is the sequence generated by the Com-DSZO algorithm for any $i\in \mathcal{V}$. Set $\mu\!=\!{1}/{\sqrt{T}}$, $\epsilon={\mu}/{2\underline{R}}$, $\psi \in (0, r^{-1}]$ and ${{\eta }_{k}}={1}/{\sqrt{k+c}}$ with $c>0$. The constant $\gamma$ is given in Appendix-A \eqref{gamma}. Then, for any optimal point $x^*\in\Omega$, we have
\begin{align}\label{thm2}
	&\sum\nolimits_{k=1}^T\mathbb{E}[f(\bar{x}_k) - f(x^*)]\leq \mathcal{O}\left(d^2\sqrt{T+c} \right) \nonumber\\
	&+\mathcal{O}\left(\frac{d^2}{\left( \psi r\omega\right) ^2}\left(\sqrt{T+c}+\ln(T+c) \right) \right).
\end{align}
\end{thm}\par 
\begin{rem}
The error bound of the Com-DSZO algorithm includes a penalty term $\mu \hat{L}$ that comes from the stochastic zeroth order (ZO) oracle, as established in Theorem 2. This reveals the non-vanishing bias of the algorithm, an inherent property of ZO methods where solution accuracy is fundamentally limited by $\mathcal{O}(\mu)$. To mitigate this, we employ an iteration-adaptive smoothing parameter $\mu = \mathcal{O}(1/\sqrt{T})$, in line with modern ZO optimization frameworks \cite{DmatchC,online1}. Crucially, although asymptotically decreasing $\mu$ (e.g., $\mu = \mathcal{O}(1/T^\alpha)$ with $\alpha > 0$) could potentially improve accuracy, such strategies require careful synchronization between the smoothing schedule and iteration complexity to maintain the $\mathcal{O}(1/\sqrt{T})$ convergence rate. This inherent accuracy-efficiency duality underscores the need to co-optimize ZO parameters in stochastic constrained settings.
\end{rem}
\subsection{Com-DSZO with stochastic batch feedback}
When implementing stochastic zeroth-order (ZO) methods in practical scenarios, the inherent variance of gradient estimation can be mitigated through mini-batch techniques \cite{b23,residual}. Specifically, we propose a \textit{doubly batched} gradient estimator with dual sample sizes $b_1$ (directional vectors) and $b_2$ (function evaluations):  
\begin{align}\label{batch-g}
	\mathbf{\tilde{g}}^{\mu}(x_k) 
	= \frac{d}{b_1b_2\mu} \sum_{\tau_1=1}^{b_1} \sum_{\tau_2=1}^{b_2} \Big[ \Big( &F(x_k+\mu u_{\tau_1}^k, \xi_{\tau_2}^k) \nonumber \\
	&- F(x_k, \xi_{\tau_2}^k) \Big) u_{\tau_1}^k \Big]
\end{align}
where $\{u_{\tau_1}^k\}_{\tau_1=1}^{b_1}$ are uniformly sampled from the unit sphere $\mathbb{S}^{d}$, and $\{\xi_{\tau_2}^k\}_{\tau_2=1}^{b_2}$ form an independent sequence replicating $\xi_k$'s distribution. This structural design yields enhanced variance control, as formalized below.
\begin{lem}[\!\!\cite{admm1}]\label{grad_var}
	Under Assumptions \ref{assum2-f}, \ref{assum3-g}, and \ref{smooth}, the gradient estimator \eqref{batch-g} satisfies:
	\begin{align}
		\mathbb{E}\left[ \big\| \mathbf{\tilde{g}}^{\mu}(x) \big\|^2 \right] \leq \frac{1}{b_1b_2} \left( 2d \big( \| \nabla f(x) \|^2 + \sigma^2 \big) + \mu^2 L^2 d^2 \right).\nonumber
	\end{align}
\end{lem}
Integrating this estimator with the Com-DSZO framework, we derive the variance-reduced variant VR-Com-DSZO. The following theorem establishes its convergence guarantees.
\begin{thm}\label{result 1 for mini-batch}
	Under the same conditions as in Theorem \ref{result 1 for smooth function}, we further suppose that $\kappa_0=1-\frac{\psi r \omega}{2}$, then, for any optimal point $x^*\in\Omega$, we have
\begin{align}\label{thm3}
	&\sum\nolimits_{k=1}^T \mathbb{E}[f(\bar{x}_k) - f(x^*)] \leq \mathcal{O}\left(\frac{d^2\sqrt{T+c}}{b_1b_2\left( \psi r \omega\right)^2 }\right) .
\end{align}
\end{thm}\par 
\begin{rem}
Theorem \ref{result 1 for mini-batch} demonstrates that the mini-batch mechanism reduces the estimator variance by a factor $1/b_1b_2$, where the batch sizes $b_1$, $b_2$ are independent of the problem dimension $d$. Remarkably, with $b_1b_2\!=\!d\sqrt{d}$, VR-Com-DSZO matches the state-of-the-art rate of centralized ZO methods \cite{b14}.
\end{rem} 
\section{Simulation}
In this section, we demonstrate and validate the theoretical results presented in this paper through numerical simulations.\par 
\textit{Optimization task:} Consider a stochastic distributed constrained optimization problem: $\min_{x \in \Omega} \frac{1}{n} \sum_{i=1}^{n} f_i(x)$,
where \( f_i(x) = \mathbb{E}_{\xi_i} \left[ \| x - \xi_i \|^2 \right] + 0.1 \| x \|_1 \) represents the nonsmooth local cost function, \( x \in \mathbb{R}^{10} \) is the global decision vector, and \( \xi_i \) denotes a random sample. For each agent \( i \), \( \xi_i \) is drawn from a Gaussian distribution, with both the mean and variance independently sampled from the interval \( [0,1] \) at each iteration. The constraint set \( \Omega \) is defined as a closed Euclidean sphere centered at the origin with radius $10$. The algorithm parameter $ \gamma=0.1$, $\psi=0.5$, $\mu=0.1$, $\epsilon=0.2$ and the step size \( \eta_k = \frac{1}{\sqrt{k + 10}} \), where \( k \) represents the iteration number. For all \( i \in \mathcal{V} \), the decision vector \( x_i \) is initialized to \( \mathbf{0}_d \). The performance of the algorithm is evaluated through the average function value: $\frac{1}{\tau} \sum_{k=1}^{\tau} f(\bar{x}_k)$,
where \( \tau \) ranges from 1 to \( T \).
\par 
\textit{Network:} We construct a distributed network consisting of $50$ computing nodes with a communication topology \(\mathcal{G} = (\mathcal{V}, \mathcal{E})\), generated using the standard network topology generation methodology from \cite{admm1}. The adjacency matrix \(\mathbf{W}\) is then derived using the Metropolis-Hastings weight rule, ensuring doubly stochasticity and spectral properties essential for distributed consensus optimization.\par 
\textit{Algorithm:} We benchmark the performance of the Com-DSZO algorithm against three prominent alternatives: (1) Distributed Stochastic Gradient Descent (DSGD) \cite{sundharram_distributed_2010}, requiring both exact first-order information and full-precision communication; (2) Distributed Stochastic Zeroth-Order (DSZO) \cite{DmatchC}, employing a two-point gradient estimator with uncompressed communication; (3) Choco-Gossip \cite{Choco-gossip}, a communication-compressed DSGD variant with error feedback that still relies on exact first-order information.\par
\begin{figure}
	\centering
	\includegraphics[width=0.8\linewidth]{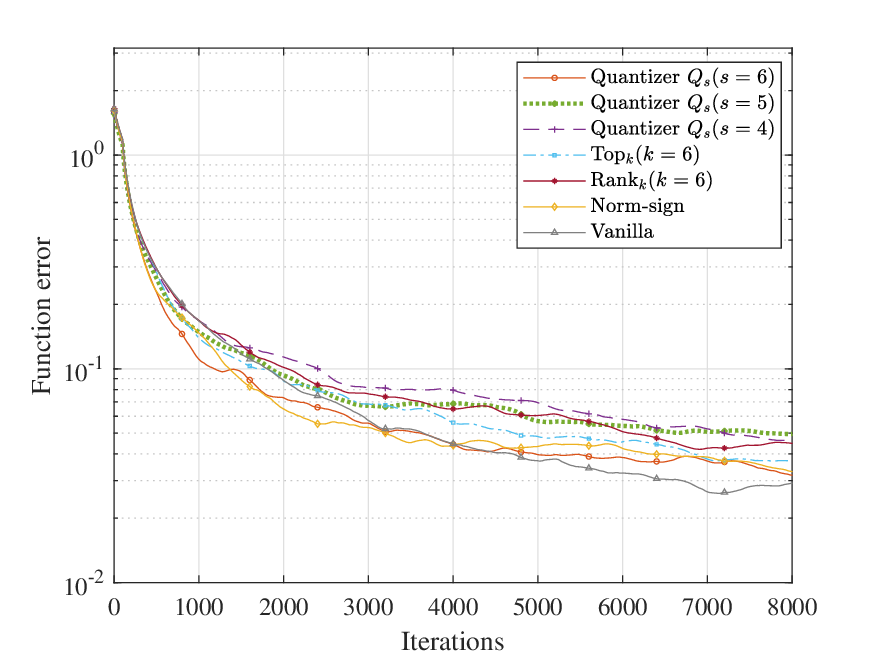}
	\caption{The effect of compression operators on algorithm performance.}
	\label{fig:tu1}
\end{figure}
\textit{Compression operators:}
We evaluate four classes of compression operators: randomized $s$-level quantizer $Q_s$ ($s=4, 5, 6$) \cite{b25} as an unbiased but contractive operator via stochastic quantization; norm-sign compressor \cite{b30} as a biased and non-contractive operator; along with two biased but contractive sparsifiers \textit{Top-$k$} and \textit{Rank-$k$} ($k=6$) \cite{onoine-choco}. \par
\textit{Result:} Our experiments lead to three primary observations:\par
\begin{itemize}
	\item As shown in Fig.\ref{fig:tu1}, the proposed Com-DSZO algorithm successfully drives the optimality gap to zero across all tested compression operators, with convergence rates closely matching the non-compressed DSZO baseline. This empirically validates the theoretical $\mathcal{O}(1/\sqrt{T})$ convergence rate established in Theorem~1. 
	\item Fig.~2 further reveals that Com-DSZO achieves performance comparable to DSGD (which utilizes exact first-order information and full-precision communication), while outperforming both DSZO and Choco-Gossip. This demonstrates Com-DSZO's robustness to gradient estimation errors and communication compression, as theoretically supported in Theorems1--2.   
	\item Additionally, Fig.~3 demonstrates that Com-DSZO reduces communication costs by approximately $49\%$--$79\%$ compared to non-compressed methods when attaining the target function error of $10^{-2}$. This substantial reduction is attributed to the error-compensated compression mechanism, which effectively minimizes data exchange overhead without compromising optimization accuracy. 
\end{itemize}
\begin{figure}
	\centering
	\includegraphics[width=0.8\linewidth]{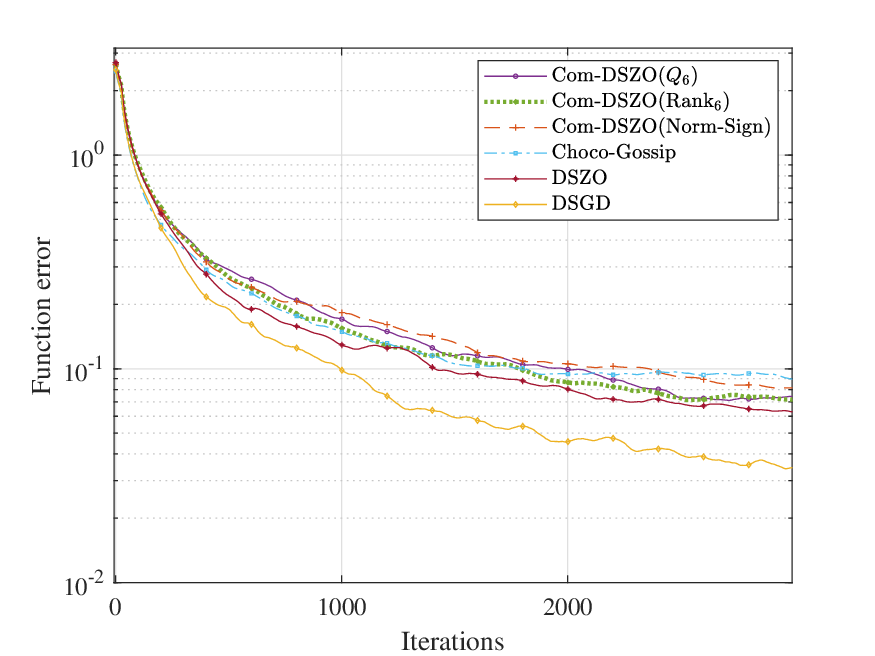}
	\caption{Comparisons of the convergence performance with Choco-SGD, DSZO, and SGD.}
	\label{fig:tu2}
\end{figure} 
\begin{figure}
	\centering
	\includegraphics[width=0.8\linewidth]{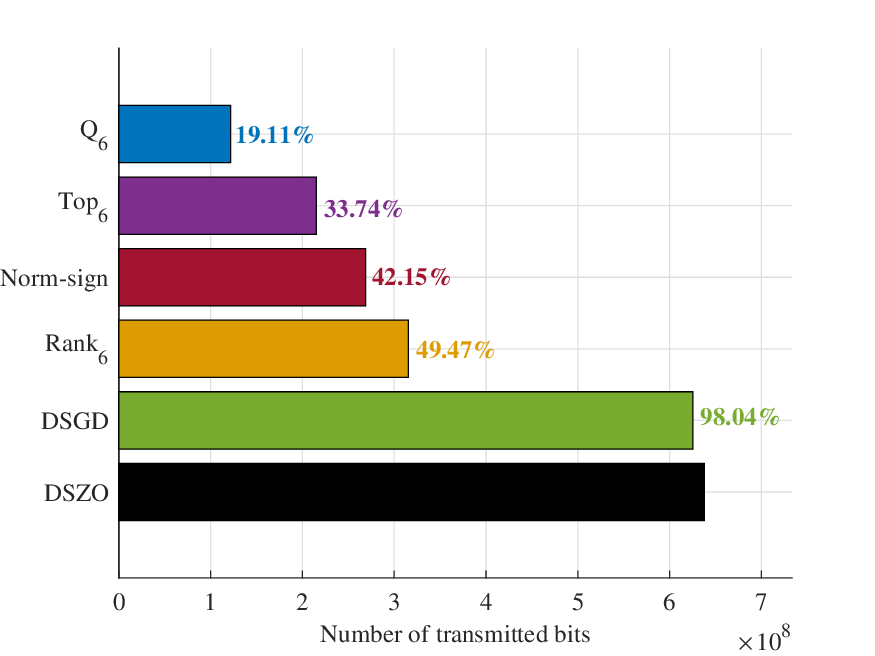}
	\caption{Comparisons of the number of transmitted bits with DSZO and SGD. }
	\label{fig:tu3}
\end{figure}

\section*{Appendix}
The following matrix inequalities are used in the proofs.\par 
\noindent\textbf{Fact 1.} For any $\mathbf{A}\in \mathbb{R}^{d\times n},\mathbf{B}\in {\mathbb{R}}^{n\times n}$,
\begin{equation}\label{Fact 1}
	\left\| \mathbf{A}\cdot \mathbf{B} \right\|_F\le \left\| \mathbf{A} \right\|_2\cdot\left\| \mathbf{B} \right\|_F.
\end{equation}\par 
\noindent\textbf{Fact 2.}  For any $\alpha >0$, $\mathbf{A}\in {{\mathbb{R}}^{d\times n}}$, $\mathbf{B}\in {{\mathbb{R}}^{d\times n}}$,
\begin{equation}\label{Fact 2}
	\left\| \mathbf{A+B} \right\|_{F}^{2}\le \left( 1+\alpha  \right){{\left\| \mathbf{A} \right\|}_{F}^2}+\left( 1+{{\alpha ^{-1}}} \right){{\left\| \mathbf{B} \right\|}_{F}^2}.
\end{equation}\par 
\subsection{Proof of Lemma \ref{smooth-e}}\label{apd1}
The proof of Lemma \ref{smooth-e} relies on the following two Lemmas. 
\begin{lem}\label{consensus error}
	Suppose Assumptions \ref{assum G}--\ref{smooth} hold. Let the sequences $\lbrace {{\mathbf{{X}}}_{k}}\rbrace _{k\geq0}$ and $\lbrace {{\mathbf{\bar{X}}}_{k}}\rbrace _{k\geq0}$ be generated by Algorithm 1. Then for any $k\geq0$, we have
		\begin{align}\label{consensus}
			& U_{1,k+1}\leq\mathbf{A}_1U_{1,k}+\mathbf{A}_2U_{2,k}+\mathbf{A}_3\eta_k^2\left\|\mathbf{G}_k^{\mu} \right\|_F^2\text{,}
		\end{align}
where $\mathbf{A}_1=\left( 1+{{\alpha }_{1}} \right)\left( 1+{{\alpha }_{2}}\right) \left(1-\gamma\delta \right) ^2$, $\mathbf{A}_2=\left( 1+{{\alpha }_{1}^{-1}} \right)\gamma^2\beta^2$ and $\mathbf{A}_3=4\left( 1+{{\alpha }_{1}} \right)\left( 1+{{\alpha }_{2}^{-1}} \right)$. Here, $\alpha_1$ and $\alpha_2$ denotes arbitrary positive numbers. 
\end{lem}\par 
\begin{IEEEproof}
	From \eqref{a} and \eqref{b} and notice that $\mathbf{\hat{X}}_0=\mathbf{B}_0=0$, it follows that $\mathbf{B}_{k+1}=\left(\mathbf{I}-\mathbf{W} \right)\mathbf{\hat{X}}_{k+1}$. Recalling \eqref{c}, we have
		\begin{align}\label{c-1}
		& U_{1,k+1}
		\leq \left\| \mathbf{\tilde{X}}_{k+1}-\mathbf{\bar{{X}}}_{k+1}\right\|_{F}^{2}\nonumber\\
		&=\left\| {\mathbf{X}_{k}}-\gamma\left(\mathbf{I}-\mathbf{W} \right)\mathbf{\hat{X}}_{k+1}-{{\eta }_{k}}\mathbf{G}_{k}^{\mu}-\mathbf{\bar{{X}}}_{k+1}\right\|_{F}^{2}\nonumber\\
		&=\Big\|\left(\left(1-\gamma \right)\mathbf{I}+\gamma \mathbf W\right) \left(  \mathbf{X}_k-\mathbf{\bar{X}}_k\right)-{{\eta }_{k}}\mathbf{G}_{k}^{\mu}\nonumber\\
		&\quad-\gamma\left(\mathbf{I}-\mathbf{W} \right)\left(\mathbf{\hat{X}}_{k+1}-\mathbf{X}_k \right)+\mathbf{\bar{X}}_k-\mathbf{\bar{X}}_{k+1}\Big\|_F^2\nonumber\\ 
		&\leq 	\left( 1+{{\alpha }_{1}^{-1}} \right)\gamma^2\beta^2 U_{2,k}\nonumber\\ 
		&\quad +\left( 1+{{\alpha }_{1}} \right)\left( 1+{{\alpha }_{2}}\right) \big\|\left(\left(1-\gamma \right)  \mathbf{I}+\gamma \mathbf W\right) \left(  \mathbf{X}_k-\mathbf{\bar{X}}_k\right)\big\|_F^2 \nonumber\\
		&\quad+\left( 1+{{\alpha }_{1}} \right)\left( 1+{{\alpha }_{2}^{-1}} \right)2\left\|\mathbf{\bar{X}}_{k+1}-\mathbf{\bar{X}}_{k} \right\|_F^2\nonumber\\
		&\quad+\left( 1+{{\alpha }_{1}} \right)\left( 1+{{\alpha }_{2}^{-1}} \right)2\eta_k^2\left\|\mathbf{G}_k^{\mu} \right\|_F^2\text{,}
	\end{align}
where the first inequality follows from the non-expansiveness of $\Pi_{\Omega_\epsilon}$, with the first equality stemming from \eqref{c} and the last inequality using \eqref{Fact 1}, \eqref{Fact 2} and $\left\|\mathbf{W}-\mathbf{I} \right\|_2=\beta$. Notice that $\left(\left(1-\gamma \right)  \mathbf{I}+\gamma \mathbf W\right) \left(  \mathbf{X}_k-\mathbf{\bar{X}}_k\right)=\left(1-\gamma \right)\mathbf{I}\left(  \mathbf{X}_k-\mathbf{\bar{X}}_k\right)+\gamma\left(  \mathbf W-\mathbf{H}\right)  \left(  \mathbf{X}_k-\mathbf{\bar{X}}_k\right)$. By $\left\| \mathbf{W} - \mathbf{H} \right\|_2 = 1 - \delta$, we obtain
\begin{align}\label{c-3}
	&\left\|\left(\left(1-\gamma \right)  \mathbf{I}+\gamma \mathbf W\right) \left(  \mathbf{X}_k-\mathbf{\bar{X}}_k\right)  \right\|_{F}^{2}\nonumber\\
	&\leq\left[\left(1-\gamma \right)\left\|\mathbf{X}_k-\mathbf{\bar{X}}_k \right\|_F+\gamma\left(1-\delta \right) \left\|\mathbf{X}_k-\mathbf{\bar{X}}_k \right\|_F  \right]^2\nonumber\\
	&\leq \left(1-\gamma\delta \right) ^2U_{1,k}\text{,}
\end{align}
where the first inequality follows from \eqref{Fact 1}. Using again the non-expansiveness of $\Pi_{\Omega_\epsilon}$ yields
\begin{align}\label{c-2}
&\left\|\mathbf{\bar{X}}_{k+1}-\mathbf{\bar{X}}_{k} \right\|_F^2\leq \left\|\mathbf{\bar{\tilde{X}}}_{k+1}-\mathbf{\bar{X}}_{k} \right\|_F^2\leq \eta_k^2\left\|\mathbf{G}_k^{\mu} \right\|_F^2\text{,}
\end{align}
where the last inequality follows from the doubly stochasticity of $\mathbf{W}$, \eqref{Fact 1} and $\left\|\mathbf{H}\right\|_2=1$. 
Substituting \eqref{c-3} and \eqref{c-2} into \eqref{c-1} gives the desired result.
\end{IEEEproof}
\begin{lem}\label{compression error}
	Suppose Assumptions \ref{assum G}--\ref{smooth} hold. Let the sequences $\lbrace {{\mathbf{{X}}}_{k}}\rbrace _{k\geq0}$ and $\lbrace {{\mathbf{\hat{X}}}_{k}}\rbrace _{k\geq0}$ be generated by Algorithm 1. Then for any $k\geq0$, we have
\begin{align}
	&\mathbb{E}_\mathcal{C}\left[U_{2,k+1}\right]\hspace{-1pt}\leq\hspace{-1pt}  \mathbf{B}_1\mathbb{E}_\mathcal{C}\left[U_{1,k}\right] \hspace{-1pt}+\hspace{-1pt}\mathbf{B}_2\mathbb{E}_\mathcal{C}\left[U_{2,k}\right]\hspace{-1pt}+\hspace{-1pt}\mathbf{B}_3\eta_k^2 \left\|\mathbf{G}_k^{\mu}\right\|_F^2\hspace{-2pt}\text{,}\nonumber
\end{align}
where $ \mathbf{B}_1=\left(1+\alpha_3^{-1} \right)\left(1+\alpha_4^{-1} \right)2\gamma^2\beta^2\left( 1-\psi r\omega\right)$, $\mathbf{B}_2=\left[ \left(1+\alpha_3 \right) +\left(1+\alpha_3^{-1} \right)\left(1+\alpha_4 \right)\gamma^2\beta^2\right] \left( 1-\psi r\omega\right)$ and $\mathbf{B}_3=\left(1+\alpha_3^{-1} \right)\left(1+\alpha_4^{-1} \right)2\left( 1-\psi r\omega\right)$.
Here, $\alpha_3$ and $\alpha_4$ are arbitrary positive numbers. 
\end{lem} \par
\begin{IEEEproof} 
	Recalling \eqref{a} and \eqref{d}, we have
	\begin{align}
		\hat{\mathbf{X}}_{k+2}-\mathbf{X}_{k+1}=	\hat{\mathbf{X}}_{k+1}+\psi\mathcal{C}\left(\mathbf{X}_{k+1}-\mathbf{\hat{X}}_{k+1}\right)-\mathbf{X}_{k+1}.\nonumber 
	\end{align}
	Denote $\mathcal{C}_r\left(\cdot \right)=\frac{1}{r}\mathcal{C}\left(\cdot \right)$. Then, from Assumption \ref{assum-c}, we can get
	\begin{align}\label{3-2}
		&\mathbb{E}_\mathcal{C}\left[ U_{2,k+1}\right]  =\mathbb{E}_\mathcal{C}\Big[\big\|( 1-\psi r) \big( \mathbf{X}_{k+1}-\hat{\mathbf{X}}_{k+1}\big)\nonumber\\
		&\quad +\psi r \big[ \big( \mathbf{X}_{k+1}-\hat{\mathbf{X}}_{k+1}\big)-\mathcal{C}_r\big(\mathbf{X}_{k+1}-\mathbf{\hat{X}}_{k+1}\big)\big] \big\| _F^2\Big]\nonumber\\
		&\leq \left( 1-\psi r\omega\right) \mathbb{E}_\mathcal{C}\Big[ \big\|\hat{\mathbf{X}}_{k+1}-\mathbf{X}_{k+1}\big\|_F^2\Big] \nonumber\\
		&\leq  \left(1+\alpha_3 \right) \left( 1-\psi r\omega\right) \mathbb{E}_\mathcal{C}\left[ U_{2,k}\right] \nonumber\\
		&\quad+ \left(1+\alpha_3^{-1} \right)\left( 1-\psi r\omega\right) \mathbb{E}_\mathcal{C}\Big[ \left\|\mathbf{X}_{k+1}-\mathbf{X}_{k}\right\|_F^2\Big]\text{,}
	\end{align}
    where the constant $\psi\in(0, r^{-1}]$, the first inequality follows from Assumption \ref{assum-c} and the convexity of F-norm, and the final one applies \eqref{Fact 2}. 
	Next, we bound the last term on the right-hand side of \eqref{3-2} as follows:
	\begin{align}\label{3-1}
		&\left\| \mathbf{{X}}_{k+1}-\mathbf{{X}}_{k} \right\|_F^2\leq \left\| \mathbf{\tilde{X}}_{k+1}-\mathbf{{X}}_{k} \right\|_F^2\nonumber\\
		&\le \left(1+\alpha_4^{-1} \right)\Big[2\gamma^2\left\|\mathbf{I}-\mathbf{W}\right\|_2^2 \left\|\mathbf{X}_k-\mathbf{\bar{X}}_k \right\|_F^2+2\eta_k^2\left\|\mathbf{G}_k^{\mu}\right\|_F^2\Big] \nonumber\\
		&\quad+\left(1+\alpha_4 \right)\gamma^2\left\|\mathbf{I}-\mathbf{W}\right\|_2^2 \left\| {\mathbf{\hat{X}}_{k+1}}-\mathbf{X}_k\right\|_F^2\nonumber\\
		&\le 2\left(1+\alpha_4^{-1} \right)\gamma^2\beta^2 \left\|\mathbf{X}_k-\mathbf{\bar{X}}_k \right\|_F^2+2\left(1+\alpha_4^{-1} \right)\eta_k^2\left\|\mathbf{G}_k^{\mu}\right\|_F^2 \nonumber\\
		&\quad+\left(1+\alpha_4 \right)\gamma^2\beta^2 \left\| {\mathbf{\hat{X}}_{k+1}}-\mathbf{X}_k\right\|_F^2\text{,}
	\end{align}
where the first inequality stems from the non-expansiveness of $\Pi_{\Omega_\epsilon}$, with the second inequality using \eqref{Fact 1}, \eqref{Fact 2} and the last inequality following from the fact that $\left\| \mathbf{I}-\mathbf{W} \right\|_2=\beta$. By combining \eqref{3-2} and \eqref{3-1}, we get the desired result.
\end{IEEEproof}
\begin{IEEEproof}[Proof of Lemma \ref{smooth-e}]
From Lemmas \ref{consensus error} and \ref{compression error}, we can obtain
	\begin{align}\label{e-diedai}
		\mathbf{e}_{k+1}&\leq \kappa_0\mathbf{e}_{k}+\mathbf{C}_3\eta_k^2\mathbb{E}\left[ \left\|\mathbf{G}_k^{\mu} \right\|_F^2\right] \text{,} 
	\end{align}
	where $\kappa_0\leq\max\left(  \mathbf{\mathbf{C}}_1,\mathbf{C}_2\right) $ with $\mathbf{C}_1=\mathbf{A}_1+\mathbf{B}_1$, $\mathbf{C}_2=\mathbf{A}_2+\mathbf{B}_2$ and $\mathbf{C}_3=\mathbf{A}_3+\mathbf{B}_3$.
	Recalling that $\left\lbrace \alpha_i\right\rbrace_{i=1,2,3,4}$ can be arbitrary positive numbers. Let $\alpha_1=\alpha_2=\frac{\gamma
		\delta}{2}$, $\alpha_3=\alpha_4^{-1}=\frac{\psi r \omega}{2}=c_0$. Then, we have $\mathbf{C}_1\leq 1-\gamma\left[\delta-\gamma\kappa_1\right]$, $\mathbf{C}_2\leq 1-c_0+\gamma\kappa_2$ and $\mathbf{C}_3\leq 4\left( 1+\frac{\gamma\delta}{2}\right)\left( 1+\frac{2}{\gamma\delta}\right)+2\left(1+c_0^{-1} \right)$=$\widetilde{\mathbf{C}}_3$, where the inequalities arise from repeated applications of the relations $\gamma^2 < \gamma < 1$ and $\left(1 + \frac{x}{2}\right)\left(1 - x\right) \leq 1 - \frac{x}{2}$. Here $\kappa_1=\frac{\delta}{4}+2(1+c_0^{-1})^2\gamma\beta^2$ and $\kappa_2= \left(1+c_0^{-1} \right)^2\gamma\beta^2+\left( 1+\frac{2}{\delta}\right)\beta^2$.
On the other hand, based on the definition of $\mathbf{G}_k^\mu$ and Lemma \ref{g}, we have 
    \begin{align}
    	\mathbb{E}\left[ \left\|\mathbf{G}_k^\mu\right\|_F^2\right]= \mathbb{E}\left[\sum\nolimits_{i=1}^n \left\|\mathbf{g}_k^\mu\right\|^2\right] \leq n d^2L_f^2.\nonumber 
    \end{align}
Let 
\begin{equation}\label{gamma}
	\gamma<\min \Big\{1,\delta\kappa_1^{-1},c_0 \kappa_2^{-1}\Big\}.
\end{equation}
It can be easily verified that $\kappa_0<1$, which further implies that $\mathbf{e}_{k+1}$ can be bounded by a perturbed contraction $\mathbf{e}_{k+1}\leq \kappa_0\mathbf{e}_{k}+\kappa_3\eta_k^{2}$ with the perturbation as $\kappa_3=
nd^2L_f^2\widetilde{\mathbf{C}}_3=\mathcal{O}(d^2c_0^{-1})$. Note that $\mathbf{e}_0=0$. The recursive update of $\mathbf{e}_{k}$ admits the formulation:
\begin{align}\label{e2}
	\mathbf{e}_{k}&\leq \kappa_3\sum\nolimits_{\tau=0}^{k-1}\kappa_0^{k-1-\tau}\eta_{\tau}^2\text{.}
\end{align}
By summing both sides of \eqref{e2} and using $\kappa_0<1$ gives
\begin{align}\label{e1}
	\sum\nolimits_{k=1}^T\mathbf{e}_{k}&\leq \frac{\kappa_3}{1-\kappa_0}\sum\nolimits_{k=0}^{T-1}\eta_{k}^2\text{.}
\end{align}\par 
Through examination of the parameters $\kappa_1$ and $\kappa_2$, it becomes evident that this relationship guarantees the satisfaction of Condition \eqref{gamma} when $\gamma\leq \frac{\delta c_0^3}{(\delta+2)c_0^2\beta^2+(c_0+1)^2\beta^2}=\mathcal{O} \left( c_0\right) $. That is, condition \eqref{gamma} is feasible. From \eqref{e-diedai}, there should be $\kappa_0\leq 1-\gamma\delta+\gamma^2\kappa_1$ and $\kappa_0\leq 1-c_0+\gamma\kappa_2$. It can be straightforward to verify that $\kappa_0=1-c_0$ is sufficient. In view of $c_0=\frac{\psi r \omega}{2}$ and $\kappa_3=
nd^2L_f^2\mathbf{\widetilde{C}}_3=\mathcal{O}(d^2c_0^{-1})$. Then, we have
\begin{align}\label{etrade-off}
	\sum\nolimits_{k=1}^T\mathbf{e}_{k}&\leq \frac{2nd^2L_f^2\mathbf{C}_3}{\psi r \omega }\sum\nolimits_{k=0}^{T-1}\eta_{k}^2\nonumber\\
	&=\mathcal{O}\left(\frac{d^2}{\left( \psi r \omega\right) ^2 }\sum\nolimits_{k=0}^{T-1}\eta_{k}^2\right) \text{,}
\end{align}
which completes the proof.
\end{IEEEproof}
\subsection{Proof of Theorem \ref{result 1 for non-smooth function}}
\begin{IEEEproof}
	Following the proof trajectory of Theorem 1, we characterize the dynamics of the optimization error evolution through the Lyapunov function $V_{k+1}(\epsilon) = \|\bar{x}_{k+1} - (1-\epsilon)x^*\|^2$
	as derived in \eqref{thm1-4}. In the nonsmooth setting, the estimation performance established in \eqref{smooth-G} becomes unattainable. Applying Lemma 1, we obtain the following upper bound for gradient estimation complexity
	\begin{align}
		\mathbb{E}\left[\big\|\mathbf{G}_{k}^{\mu}\big\|^2\right] 
		&= \sum\nolimits_{i=1}^n \mathbb{E}_{\xi_k}\left[\mathbb{E}_{u_k}\left[\big\|\mathbf{g}_{i,k}^{\mu}\big\|^2\right]\right]\leq nd^2L_f^2.
	\end{align}
	Combining this bound with \eqref{thm1-4} yields
	\begin{align}\label{thm2-2}
		\mathbb{E}[V_{k+1}(\epsilon)] &\leq \mathbb{E}[V_k(\epsilon)] - 2\eta_k\mathbb{E}[f(\bar{x}_k) - f(x^*)] \nonumber\\
		&\quad + \left(\frac{10}{n} + 5\right)d^2L_f^2\eta_k^2 + \left(\frac{2}{n^2} + \frac{1}{n}\right)4\gamma^2\beta^2\mathbf{e}_k \nonumber\\
		&\quad + (\bar{R}\epsilon + 2\mu)2L_f\eta_k + \frac{2L_m\eta_k}{n}\mathbb{E}[U_{1,k}].
	\end{align}
	Based on $\mathbb{E}[U_{1,k}] < \mathbf{e}_k$ and Lemma \ref{smooth-e}, we deduce
	\begin{align}\label{thm2-4}
		&\sum\nolimits_{k=1}^T\mathbb{E}[f(\bar{x}_k) - f(x^*)]\leq \frac{2\bar{R}^2}{\eta_T} + (\bar{R}\epsilon + 2\mu)L_fT \nonumber\\
		&+ \left(\frac{5}{n} + \frac{5}{2}\right)d^2L_f^2\sum\nolimits_{k=1}^T\eta_k + \frac{L_m}{n}\frac{\kappa_3}{1-\kappa_0}\sum\nolimits_{k=1}^T\mathbf{e}_k \nonumber\\
		& + \left(\frac{1}{n^2} + \frac{1}{2n}\right)4\gamma^2\beta^2\sum\nolimits_{k=1}^T\frac{\mathbf{e}_k}{\eta_k}.
	\end{align}
	Final convergence rate is obtained by combining \eqref{thm1-6} and \eqref{thm2-4} with parameter choices $\mu = 1/\sqrt{T}$ and $\epsilon = \mu/(2\underline{R})$:
	\begin{align}\label{thm2-5}
		&\sum\nolimits_{k=1}^T\mathbb{E}[f(\bar{x}_k) - f(x^*)]\leq 2\bar{R}^2\sqrt{T+c} +\left( 2+{\frac{\overline{R}}{2\underline{R}}}\right)L_f\sqrt{T} \nonumber\\
		&+ \left(\frac{5}{n} + \frac{5}{2}\right)d^2L_f^2\sqrt{T+c} + \frac{L_m\mathbf{C}_3L_f^2d^2}{n(1-\kappa_0)}\ln(T+c) \nonumber\\
		&+ \left(\frac{2}{n^2} + \frac{1}{n}\right)\frac{4\gamma^2\beta^2\mathbf{C}_3L_f^2d^2}{1-\kappa_0}\sqrt{T+c}.
	\end{align}
Let $\kappa_0=1-c_0=1-\frac{\psi r \omega}{2}$. Similar to the way to get \eqref{corollary1}, we have 
\begin{align}\label{thm2-6}
&\sum\nolimits_{k=1}^T\mathbb{E}[f(\bar{x}_k) - f(x^*)]\leq \mathcal{O}\left(d^2\sqrt{T+c} \right) \nonumber\\
&+\mathcal{O}\left(\frac{d^2}{\left( \psi r\omega\right) ^2}\left(\sqrt{T+c}+\ln(T+c) \right) \right).
\end{align}
\end{IEEEproof}
\subsection{Proof of Theorem 3}
\begin{IEEEproof}
	Following the methodology used in proving \eqref{thm1-5} and relying on Lemma \ref{g}, we establish
	\begin{align}\label{thm3-1}
		&\mathbb{E}\left[V_{k+1}(\epsilon)\right]\leq \mathbb{E}\left[V_k(\epsilon)\right] - 2\eta_k\big[1 - \widetilde{M}_1\eta_k\big] \mathbb{E}\left[f(\bar{x}_k) - f(x^*)\right] \nonumber\\
		&+\left(\frac{2}{n^2} + \frac{1}{n}\right)4\gamma^2\beta^2\mathbf{e}_k + \left(\frac{2L_m\eta_k}{n} + 4dL_m^2\widetilde{M}_0\eta_k^2\right)\mathbb{E}[U_{1,k}] \nonumber\\
		&+(\bar{R}\epsilon + 2\mu)2L_f\eta_k + (4nd\hat{\sigma}^2 + nd^2L_m^2\mu^2)\frac{\widetilde{M}_0\eta_k^2}{2},
	\end{align}
	where $\widetilde{M}_0 = \frac{M_0}{b_1b_2}$ and $\widetilde{M}_1 = 4ndL_m\widetilde{M}_0$. Choose the step-size $\eta_k$ as $1/\big( \widetilde{M}_1(k+c)^a\big) $ with $c>0$. Rearranging the terms in \eqref{thm3-1} leads to
	\begin{align}\label{thm3-2}
		& \mathbb{E}\left[f(\bar{x}_k) - f(x^*)\right]
		\leq \frac{\widetilde{M}_1(k+c)^{2a}}{2((k+c)^{a}-1)}\big\{\mathbb{E}\left[V_k(\epsilon)\right] -\mathbb{E}\left[V_{k+1}(\epsilon)\right]\big\}\nonumber\\
		&+\left(\frac{2}{n^2} + \frac{1}{n}\right)2\gamma^2\beta^2\widetilde{M}_1\mathbf{e}_k\cdot\frac{(k+c)^{2a}}{(k+c)^{a}-1} \nonumber\\
		&+ \frac{L_m}{n}\mathbb{E}[U_{1,k}] +\left( \frac{L_m}{n}+ \frac{2dL_m^2\widetilde{M}_0}{\widetilde{M}_1}\right)  \cdot\frac{\mathbb{E}[U_{1,k}]}{(k+c)^{a}-1}\nonumber\\
		&+(\bar{R}\epsilon + 2\mu)L_f \cdot\frac{(k+c)^{a}}{(k+c)^{a}-1} \nonumber\\
		&+ (4nd\hat{\sigma}^2+ nd^2L_m^2\mu^2)\frac{\widetilde{M}_0}{4\widetilde{M}_1}\cdot\frac{1}{(k+c)^{a}-1}
	\end{align}
	Applying Lemma \ref{grad_var} through analogous derivations to \eqref{e2}, we derive
	\begin{align}\label{thm3-3}
		\sum\nolimits_{k=1}^T \mathbf{e}_k \leq \frac{\tilde{\kappa}_3}{1 - \kappa_0}\sum\nolimits_{k=0}^{T-1}\eta_k^2,
	\end{align}
	where $\tilde{\kappa}_3 = \frac{\kappa_3}{b_1b_2}$. Observe that $\frac{(k+c)^{2a}}{(k+c)^a-1} \leq \big(2+(k+c)^a\big)$, $\sum_{k=0}^T \frac{1}{(k + c)^a - 1} < 2(T + c)^{1 - a}$ and $\underline{R}$ denotes the lower bound of parameter $R$. Set $a=1/2$, $\mu=1/{\sqrt{T}}$ and $\epsilon={\mu}/{2\underline{R}}$. Recalling \eqref{1-2} and summing over $k = 1$ to $T$, we deduce:
	\begin{align}\label{thm3-4}
		&\sum\nolimits_{k=1}^T \mathbb{E}[f(\bar{x}_k) - f(x^*)] \leq 2\bar{R}^2\widetilde{M}_1\left(2 + \sqrt{T + c}\right) \nonumber\\
		&+ \left(\frac{2}{n^2} + \frac{1}{n}\right)4\gamma^2\beta^2\widetilde{M}_1T + (\bar{R}\underline{R} + 1)2L_fT \nonumber\\
		&+ \left(\frac{2}{n^2} + \frac{1}{n}\right)\frac{2\gamma^2\beta^2\tilde{\kappa}_3\sqrt{T+c}}{1-\kappa_0 } + \left( 2+{\frac{\overline{R}}{2\underline{R}}}\right)L_f\sqrt{T+c} \nonumber\\
		& +\frac{{\kappa}_3\sqrt{T+c}}{2n^2dM_0\left( 1-\kappa_0\right) }+\frac{L_m\tilde{\kappa}_3\sqrt{T+c}}{nM_1\left( 1-\kappa_0\right) }+\frac{L_m\tilde{\kappa}_3\ln(T+c)}{n\widetilde{M_1}^2\left( 1-\kappa_0\right) } \nonumber\\
		&+\frac{\hat{\sigma}^2\sqrt{T+c}}{2L_m}+\frac{dL_m}{4\sqrt{T}}.
	\end{align}
	Given that $\tilde{\kappa_3}\propto d^2/b_1b_2c_0$, $M_1\propto d$ and $\kappa_0=1-c_0$. We can develop the desired result in \eqref{thm3}.
\end{IEEEproof}
\bibliographystyle{IEEEtran}
\bibliography{ref}

\begin{thebibliography}{10}
\providecommand{\url}[1]{#1}
\csname url@samestyle\endcsname
\providecommand{\newblock}{\relax}
\providecommand{\bibinfo}[2]{#2}
\providecommand{\BIBentrySTDinterwordspacing}{\spaceskip=0pt\relax}
\providecommand{\BIBentryALTinterwordstretchfactor}{4}
\providecommand{\BIBentryALTinterwordspacing}{\spaceskip=\fontdimen2\font plus
\BIBentryALTinterwordstretchfactor\fontdimen3\font minus
  \fontdimen4\font\relax}
\providecommand{\BIBforeignlanguage}[2]{{%
\expandafter\ifx\csname l@#1\endcsname\relax
\typeout{** WARNING: IEEEtran.bst: No hyphenation pattern has been}%
\typeout{** loaded for the language `#1'. Using the pattern for}%
\typeout{** the default language instead.}%
\else
\language=\csname l@#1\endcsname
\fi
#2}}
\providecommand{\BIBdecl}{\relax}
\BIBdecl

\bibitem{smartgrid}
Y.~Wang, S.~Liu, B.~Sun, and X.~Li, ``A distributed proximal primal–dual
  algorithm for energy management with transmission losses in smart grid,''
  \emph{IEEE Trans. Ind. Informat.}, vol.~18, no.~11, pp. 7608--7618, 2022.

\bibitem{machinelearning}
A.~Nedi{\'c}, A.~Olshevsky, and C.~A. Uribe, ``Fast convergence rates for
  distributed non-bayesian learning,'' \emph{IEEE Trans. Autom. Control},
  vol.~62, no.~11, pp. 5538--5553, 2017.

\bibitem{gt}
J.~Lei, P.~Yi, J.~Chen, and Y.~Hong, ``Distributed variable sample-size
  stochastic optimization with fixed step-sizes,'' \emph{IEEE Trans. Autom.
  Control}, vol.~67, no.~10, pp. 5630--5637, 2022.

\bibitem{sundharram_distributed_2010}
S.~S. Ram, A.~Nedi{\'c}, and V.~V. Veeravalli, ``Distributed stochastic
  subgradient projection algorithms for convex optimization,'' \emph{J. Optim.
  Theory Appl.}, vol. 147, no.~3, pp. 516--545, 2010.

\bibitem{DSGT}
S.~Pu and A.~Nedi{\'c}, ``A distributed stochastic gradient tracking method,''
  in \emph{Proc. IEEE Conf. Decis. Control}, 2018, pp. 963--968.

\bibitem{DPGD}
S.~Liu, Z.~Qiu, and L.~Xie, ``Convergence rate analysis of distributed
  optimization with projected subgradient algorithm,'' \emph{Automatica},
  vol.~83, pp. 162--169, 2017.

\bibitem{b4}
H.~Li, Q.~Lü, and T.~Huang, ``Distributed projection subgradient algorithm
  over time-varying general unbalanced directed graphs,'' \emph{IEEE Trans.
  Autom. Control}, vol.~64, no.~3, pp. 1309--1316, 2019.

\bibitem{ADMM}
W.~Shi, Q.~Ling, K.~Yuan, G.~Wu, and W.~Yin, ``On the linear convergence of the
  admm in decentralized consensus optimization,'' \emph{IEEE Trans. Signal
  Process.}, vol.~62, no.~7, pp. 1750--1761, 2014.

\bibitem{b38}
L.~Romao, K.~Margellos, G.~Notarstefano, and A.~Papachristodoulou,
  ``Subgradient averaging for multi-agent optimisation with different
  constraint sets,'' \emph{Automatica}, vol. 131, p. 109738, 2021.

\bibitem{b37}
S.~Cheng, S.~Liang, Y.~Fan, and Y.~Hong, ``Distributed gradient tracking for
  unbalanced optimization with different constraint sets,'' \emph{IEEE Trans.
  Autom. Control}, vol.~68, no.~6, pp. 3633--3640, 2023.

\bibitem{flexman-onepoint}
A.~D. Flaxman, A.~T. Kalai, and H.~B. McMahan, ``Online convex optimization in
  the bandit setting: gradient descent without a gradient,'' in \emph{Proc.
  16th Annu. ACM-SIAM Symp. Discrete Algorithms}, 2005, pp. 385--394.

\bibitem{nesterov2017}
Y.~Nesterov and V.~Spokoiny, ``Random gradient-free minimization of convex
  functions,'' \emph{Found. Comput. Math.}, vol.~17, no.~2, pp. 527--566, 2017.

\bibitem{b14}
O.~Shamir, ``An optimal algorithm for bandit and zero-order convex optimization
  with two-point feedback,'' \emph{J. Mach. Learn. Res.}, vol.~18, no.~1, pp.
  1703--1713, 2017.

\bibitem{duchi}
J.~C. Duchi, M.~I. Jordan, M.~J. Wainwright, and A.~Wibisono, ``Optimal rates
  for zero-order convex optimization: The power of two function evaluations,''
  \emph{IEEE Trans. Inform. Theory}, vol.~61, no.~5, pp. 2788--2806, 2015.

\bibitem{DmatchC}
D.~Yuan, L.~Wang, A.~Proutiere, and G.~Shi, ``Distributed zeroth-order
  optimization: Convergence rates that match centralized counterpart,''
  \emph{Automatica}, vol. 159, p. 111328, 2024.

\bibitem{nonsmooth}
Y.~Wang, W.~Zhao, Y.~Hong, and M.~Zamani, ``Distributed subgradient-free
  stochastic optimization algorithm for nonsmooth convex functions over
  time-varying networks,'' \emph{SIAM J. Control Optim.}, vol.~57, no.~4, pp.
  2821--2842, 2019.

\bibitem{noncvx-dso-zo}
X.~Yi, S.~Zhang, T.~Yang, and K.~H. Johansson, ``Zeroth-order algorithms for
  stochastic distributed nonconvex optimization,'' \emph{Automatica}, vol. 142,
  p. 110353, 2022.

\bibitem{online1}
Y.~Pang and G.~Hu, ``Randomized gradient-free distributed online optimization
  via a dynamic regret analysis,'' \emph{IEEE Trans. Autom. Control}, vol.~68,
  no.~11, pp. 6781--6788, 2023.

\bibitem{online2}
X.~Yi, X.~Li, T.~Yang, L.~Xie, T.~Chai, and K.~H. Johansson, ``Distributed
  bandit online convex optimization with time-varying coupled inequality
  constraints,'' \emph{IEEE Trans. Autom. Control}, vol.~66, no.~10, pp.
  4620--4635, 2021.

\bibitem{online3}
D.~Yuan, Y.~Hong, D.~W.~C. Ho, and S.~Xu, ``Distributed mirror descent for
  online composite optimization,'' \emph{IEEE Trans. Autom. Control}, vol.~66,
  no.~2, pp. 714--729, 2021.

\bibitem{b22}
X.~Yi, S.~Zhang, T.~Yang, T.~Chai, and K.~H. Johansson, ``Linear convergence of
  first- and zeroth-order primal–dual algorithms for distributed nonconvex
  optimization,'' \emph{IEEE Trans. Autom. Control}, vol.~67, no.~8, pp.
  4194--4201, 2022.

\bibitem{b23}
D.~Hajinezhad, M.~Hong, and A.~Garcia, ``Zone: Zeroth-order nonconvex
  multiagent optimization over networks,'' \emph{IEEE Trans. Autom. Control},
  vol.~64, no.~10, pp. 3995--4010, 2019.

\bibitem{zero-ordernonconvex}
Y.~Tang, J.~Zhang, and N.~Li, ``Distributed zero-order algorithms for nonconvex
  multiagent optimization,'' \emph{IEEE Trans. Control Netw. Syst.}, vol.~8,
  no.~1, pp. 269--281, 2021.

\bibitem{b25}
D.~Alistarh, T.~Hoefler, M.~Johansson, N.~Konstantinov, S.~Khirirat, and
  C.~Renggli, ``The convergence of sparsified gradient methods,'' in \emph{Adv.
  Neural Inf. Process. Syst.}, vol.~31, 2018, pp. 5977--5987.

\bibitem{b26}
D.~Alistarh, D.~Grubic, J.~Z. Li, R.~Tomioka, and M.~Vojnovic, ``Qsgd:
  communication-efficient sgd via gradient quantization and encoding,'' in
  \emph{Proc. 31st Int. Conf. Neural Inf. Process. Syst.}, 2017, pp.
  1707--1718.

\bibitem{unbiased1}
D.~Yuan, B.~Zhang, D.~W.~C. Ho, W.~X. Zheng, and S.~Xu, ``Distributed online
  bandit optimization under random quantization,'' \emph{Automatica}, vol. 146,
  p. 110590, 2022.

\bibitem{Choco-gossip}
A.~Koloskova, S.~Stich, and M.~Jaggi, ``Decentralized stochastic optimization
  and gossip algorithms with compressed communication,'' in \emph{Proc. 36th
  Int. Conf. Mach. Learn.}, vol.~97, 2019, pp. 3478--3487.

\bibitem{onoine-choco}
X.~Cao and T.~Ba\c{s}ar, ``Decentralized online convex optimization with
  compressed communications,'' \emph{Automatica}, vol. 156, p. 111186, 2023.

\bibitem{b29}
Y.~Liao, Z.~Li, K.~Huang, and S.~Pu, ``A compressed gradient tracking method
  for decentralized optimization with linear convergence,'' \emph{IEEE Trans.
  Autom. Control}, vol.~67, no.~10, pp. 5622--5629, 2022.

\bibitem{b30}
J.~Zhang, K.~You, and L.~Xie, ``Innovation compression for
  communication-efficient distributed optimization with linear convergence,''
  \emph{IEEE Trans. Autom. Control}, vol.~68, no.~11, pp. 6899--6906, 2023.

\bibitem{b33}
X.~Yi, S.~Zhang, T.~Yang, T.~Chai, and K.~H. Johansson, ``Communication
  compression for distributed nonconvex optimization,'' \emph{IEEE Trans.
  Autom. Control}, vol.~68, no.~9, pp. 5477--5492, 2023.

\bibitem{admm1}
X.~Gao, B.~Jiang, and S.~Zhang, ``On the information-adaptive variants of the
  admm: An iteration complexity perspective,'' \emph{J. Sci. Comput.}, vol.~76,
  no.~1, pp. 327--363, 2018.

\bibitem{nonconvexcom}
J.~Li, C.~Li, J.~Fan, and T.~Huang, ``Online distributed stochastic gradient
  algorithm for nonconvex optimization with compressed communication,''
  \emph{IEEE Trans. Autom. Control}, vol.~69, no.~2, pp. 936--951, 2024.

\bibitem{residual}
Y.~Zhang, Y.~Zhou, K.~Ji, and M.~M. Zavlanos, ``A new one-point
  residual-feedback oracle for black-box learning and control,''
  \emph{Automatica}, vol. 136, p. 110006, 2022.

\end{thebibliography}
\end{document}